\documentclass[journal]{new-aiaa}
\usepackage[utf8]{inputenc}

\usepackage{bm}
\usepackage{amsmath}

\usepackage{amsthm}
\usepackage{color,soul}
\usepackage{todonotes}
\usepackage{colortbl}
\definecolor{myblue}{rgb}{0, 0.23, 0.64}

\usepackage{hyperref}
\hypersetup{
    colorlinks=true,
    linkcolor=myblue,
    filecolor=magenta,
    citecolor=myblue,
    urlcolor=myblue,
}

\usepackage{graphicx}
\usepackage[version=4]{mhchem}
\usepackage[group-separator={,}]{siunitx}
\usepackage{longtable}
\usepackage{gensymb}
\usepackage{mathtools}
\usepackage{float}
\usepackage{subcaption}
\usepackage{mwe}
\usepackage{enumitem}
\usepackage{multirow}
\usepackage{booktabs}

\usepackage{algorithm}
\usepackage{algpseudocode}

\usepackage{lscape}
\usepackage{adjustbox}
\usepackage{dsfont}
\usepackage{xspace}
\usepackage{tcolorbox}
\usepackage[flushleft]{threeparttable}
\usepackage{mdframed}
\mdfsetup{
middlelinecolor=red,
middlelinewidth=2pt,
backgroundcolor=red!5}

\newtheorem{theorem}{Theorem}

\theoremstyle{definition}

\newtheorem{example}{Illustrative Example}

\theoremstyle{remark}
\newtheorem{remark}{Remark}

\let\svthefootnote\thefootnote
\newcommand\freefootnote[1]{
  \let\thefootnote\relax
  \footnotetext{#1}
  \let\thefootnote\svthefootnote
}

\title{Optimal Satellite Constellation Configuration Design: A Collection of Mixed Integer Linear Programs}

\author{David O. Williams Rogers\footnote{Ph.D. Student, Department of Mechanical, Materials and Aerospace Engineering, Student Member AIAA.}}
\affil{West Virginia University, Morgantown, West Virginia, 26506}
\author{Dongshik Won\footnote{Director, Future Innovation Research Team, TelePIX; Ph.D. Candidate, Department of Aerospace Engineering, KAIST, Member AIAA.}}
\affil{TelePIX, Seoul, 07330, Republic of Korea,\\
Korea Advanced Institute of Science and Technology, Daejeon 34141, Republic of Korea}
\author{Dongwook Koh\footnote{Director, Satellite Development Team.}, and Kyungwoo Hong\footnote{Principal Research Engineer, Future Innovation Research Team; now affiliated with the Agency for Defense Development, Republic of Korea.}}
\affil{TelePIX, Seoul, 07330, Republic of Korea}
\author{Hang Woon Lee\footnote{Assistant Professor, Department of Mechanical, Materials and Aerospace Engineering; hangwoon.lee@mail.wvu.edu. Member AIAA (Corresponding Author).}}
\affil{West Virginia University, Morgantown, West Virginia, 26506}

\begin{document}

\freefootnote{This paper is a substantially revised and expanded version of the Paper AIAA 2023-4658, presented at the 2023 ASCEND, Las Vegas, NV, October 23-25, 2023. It offers new results and a better description of the materials.}

\maketitle

\begin{abstract} 
Designing satellite constellation systems involves complex multidisciplinary optimization in which coverage serves as a primary driver of overall system cost and performance. Among the various design considerations, constellation configuration, which dictates how satellites are placed and distributed in space relative to each other, predominantly determines the resulting coverage. In constellation configuration design, coverage may be treated either as an optimization objective or as a constraint, depending on mission goals. State-of-the-art literature addresses each mission scenario on a case-by-case basis, employing distinct assumptions, modeling techniques, and solution methods. While such problem-specific approaches yield valuable insights, users often face implementation challenges when performing trade-off studies across different mission scenarios, as each scenario must be handled distinctly. In this paper, we propose a collection of five mixed-integer linear programs that are of practical significance, extensible to more complex mission narratives through additional constraints, and capable of obtaining provably optimal constellation configurations. The framework can handle various metrics and mission scenarios, such as percent coverage, average or maximum revisit times, a fixed number of satellites, spatiotemporally varying coverage requirements, and static or dynamic targets. The paper presents several case studies and comparative analyses to demonstrate the versatility of the proposed framework.
\end{abstract}

\section*{Nomenclature}
{\renewcommand\arraystretch{1.0}
\noindent\begin{longtable*}{@{}l @{\quad=\quad} l@{}}
\multicolumn{2}{@{}l}{Sets}\\
$\mathcal{J}$ & Set of orbital slots\\
$\mathcal{P}$ & Set of targets\\
$\mathcal{T}$ & Set of time steps\\
\multicolumn{2}{@{}l}{Subscripts and indexing}\\
$j,k$ & Orbital slot index\\
$p$ & Target index\\
$t$ & Time step index \\
\multicolumn{2}{@{}l}{Parameters}\\
$b$ & Coverage timeline\\
$c$ & Orbital slot cost\\
$D$ & Spatiotemporal coverage parameter\\
$N$ & Number of satellites\\
$r$ & Coverage threshold parameter\\
$V$ & Boolean visibility parameter \\
$W$ & Boolean inter-satellite link parameter\\
$z$ & Upper bound maximum revisit time parameter\\
$\gamma$ & Average revisit time upper bound parameter\\
$\pi$ & Observation reward \\
\multicolumn{2}{@{}l}{Decision variables}\\
$s$ & Average revisit time auxiliary variable\\
$w$ & Coverage gap duration indicator variable\\
$x$ & Satellite location decision variable \\
$y$ & Coverage state decision variable \\
$Z$ & Maximum revisit time decision variable\\
$\alpha$ & Average revisit time upper bound decision variable\\
$\ell$ & Coverage gap start indicator variable\\
\end{longtable*}}

\section{Introduction}

Satellite constellations constitute a critical component of spaceborne infrastructure, supporting a wide range of applications in scientific research and civil services, including remote sensing \cite{stephens_bams_2002,paek_asr_2018,xue_ijrs_2008}, telecommunications \cite{leyva_ieee_2020,kodheli_ieee_2017,evans_ieee_2014}, and positioning, navigation, and timing services \cite{gps,glonass,galileo}. Beyond these established capabilities, emerging applications of constellation-based systems include cislunar space situational awareness (SSA) \cite{clareson_jsr_2025,patel_jas_2024}, orbital debris remediation \cite{williams_asr_2025}, Internet of Things connectivity \cite{qu_access_2017}, and space-based solar power transmission to terrestrial power grids \cite{mankins_aa_2002,li_aa_2017}. As distributed systems, satellite constellations offer inherent advantages over traditional monolithic architectures, including enhanced resilience to individual satellite failures, improved spatiotemporal resolution, and increased overall reliability and mission performance. Furthermore, the increased availability of commercial rideshare launch opportunities and streamlined satellite manufacturing have substantially lowered barriers to space access. Collectively, these developments have strengthened the case for constellation-based architectures, expanding the capabilities of distributed space systems and establishing a foundation for more scalable and cost-effective missions.

Satellite constellation system design entails a complex multidisciplinary optimization process encompassing technical, social, and policy-related dimensions (\textit{e.g.}, spacecraft design, constellation configuration, ground segment, frequency spectrum utilization, decommissioning procedures). The tightly coupled nature of each discipline demands an integrated system design approach to address the interdependencies between system components. In response, the literature often adopts multidisciplinary design optimization (MDO). MDO is a methodology tailored for optimizing systems composed of multiple interacting disciplines, each with its own local design variables and coupling variables, facilitating information exchange across disciplines and capturing their interrelationships \cite{martins_aiaa_2013}. A system-level optimizer coordinates the individual disciplines, seeking to optimize a global objective while ensuring both system-wide (global) and discipline-specific (local) feasibility constraints are satisfied. Notable applications of MDO for constellation system design optimization include minimizing the total cost by tackling the constellation configuration, spacecraft design, and launch manifest \cite{budianto_jsr_2004}, and minimizing the total system mass cost considering constellation configuration and spacecraft design \cite{shi_smo_2018}.
    
Critical factors in the constellation system design process, such as mission performance and cost, are significantly influenced by coverage. For Earth observation and surveillance missions, coverage dictates the amount of data that can be collected from targets, whereas in telecommunications, it determines the number of users served and the quality of service provided. For each mission context, coverage requirements are jointly characterized by two dimensions: the geographical region of interest (regional or global) and the desired temporal resolution (continuous or discontinuous). For discontinuous coverage specifically, Ref.~\cite{wertz_constellation} identifies key metrics for evaluating mission performance, including \textit{percentage coverage}---the percentage of time, with respect to the mission horizon, that the constellation covers the target---\textit{maximum revisit time} (MRT)---the longest coverage gap between successive observations---and \textit{average revisit time} (ART)---the mean coverage gap between observations. The constellation's coverage depends on various factors, for instance, the payload characteristics (\textit{e.g.,} remote sensing sensors or communications' antennas), the number of satellites and their orbital parameters, and exogenous factors (\textit{e.g.,} cloud coverage). Among the various factors involved, the \textit{constellation configuration}, which dictates the number of satellites, their geometric arrangement, and their orbital characteristics, plays a central role in determining the coverage delivered by the system \cite{shi_smo_2018,budianto_jsr_2004}. 

The extensive research conducted in the literature, as introduced in Sec.~\ref{sec:lit_rev}, proves the value of constellation configuration design using coverage as a key figure or merit; however, the field exhibits considerable fragmentation that remains inadequately addressed. In particular, current research treats each mission scenario through heterogeneous modeling frameworks and computational approaches ranging from nonlinear geometrical models to integer linear programs (ILP) that impede direct comparison and adaptation across disparate mission contexts. While individual studies achieve rigor within their respective scopes, mission designers encounter significant difficulties when evaluating design alternatives, as each problem instance necessitates independent formulation and solution methodology. These isolated efforts are frequently constrained by assumptions regarding orbital geometry, satellite distribution requirements, or target characteristics, which collectively create barriers to systematic design exploration and comparative mission assessment. Consequently, a unified and generalizable framework for constellation configuration optimization across diverse coverage objectives and mission specifications remains absent from the literature.

In response to this gap, we propose a collection of five mixed-ILP (MILP) formulations, each tailored to a specific coverage scenario. All formulations are MILP at large, admitting provably optimal solutions via exact methods (\textit{e.g.,} branch and bound) or commercial solvers (\textit{e.g.,} Gurobi Optimizer, CPLEX). They accommodate complex spatiotemporal coverage requirements for single or multiple, static or dynamic targets, and impose no prescribed (a)symmetry constraints on satellite distributions. The first, the Set Covering Location Problem (SCLP; introduced in Ref.~\cite{lee2020satellite}), determines the minimum-cost constellation configuration that ensures continuous coverage over the targets. The second, the Partial Set Covering Location Problem (PSCLP), determines the minimum-cost constellation configuration that achieves at least a specified coverage threshold over the targets. The third, the Maximal Covering Location Problem (MCLP; introduced in Ref.~\cite{lee2023mclp}), determines, for a given number of satellites, the constellation configuration that maximizes the target's observation rewards. The fourth, the Minimal Maximum Revisit Time Problem (MMRT), determines, for a given number of satellites, the constellation configuration that minimizes the MRT over the targets. The fifth, the Minimal Average Revisit Time Problem (MART), determines, for a given number of satellites, the constellation configuration that minimizes the ART over the targets.

\begin{figure}[h!]
    \centering
    \includegraphics[width=0.95\linewidth]{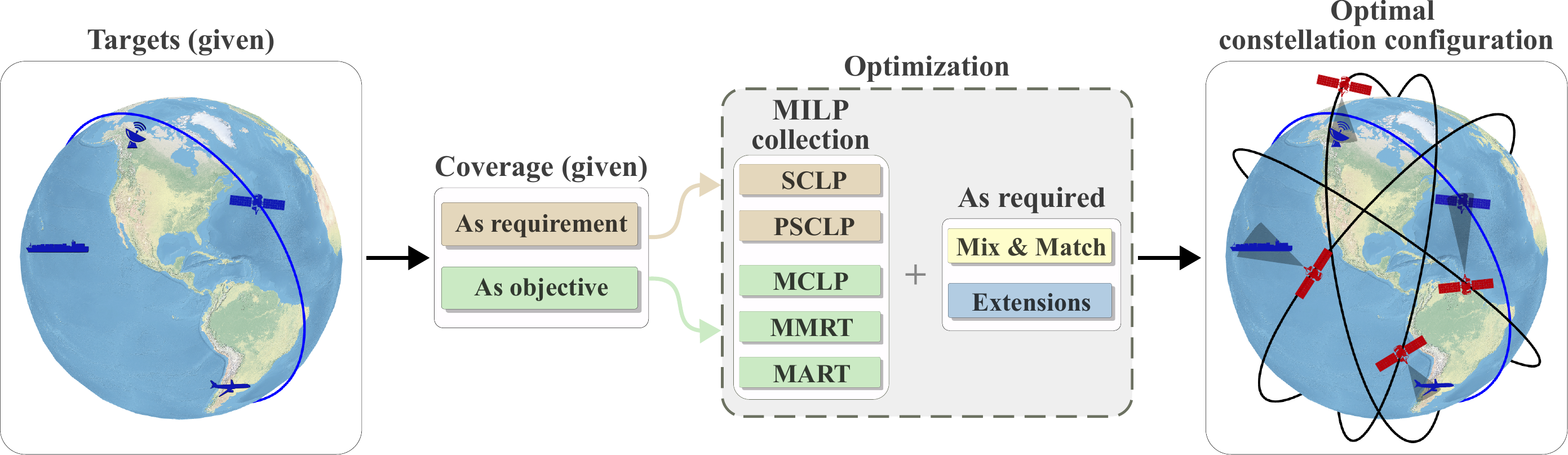}
    \caption{Satellite constellation configuration design optimization flow. The proposed MILP formulations are categorized by whether coverage is treated as a requirement or an objective.}
    \label{fig:overarching}
\end{figure}

We envision this paper to make the constellation configuration design optimization more coherent and accessible by providing mission designers with clear guidance for easy and rapid preliminary assessment of a given mission scenario (\textit{e.g.,} remote sensing, telecommunications, SSA). Given a set of targets, the mission designer specifies whether coverage is treated as a requirement (constraint) or an objective and selects the corresponding formulation. Solving the chosen problem yields the optimal constellation configuration that satisfies the stated coverage specifications over the targets. Recognizing that the five formulations may not comprehensively address all mission contexts and design scenarios, we demonstrate how mission designers can leverage MILP modeling techniques to perform mix-and-match operations between objective functions and constraints, as well as develop customized extensions suited to their specific requirements (see Fig.~\ref{fig:overarching}). For example, we derive formulations that minimize constellation configuration cost subject to user-defined upper bounds on maximum or average revisit time. We also extend the applicability of the formulations through mission parameters that account for dynamic targets and satellite deployment costs, and through additional constraints that ensure robust inter-satellite link (ISL) topologies. It is noteworthy that SCLP, PSCLP, and MCLP correspond directly to classical facility location problems (FLPs), which are well-established optimization problems in the operations research literature. These problems determine optimal facility locations while accounting for constraints such as location and transportation costs, and customer demands. For more detailed information on formulation variants and solution approaches, we refer readers to the FLP literature. This paper extends a preliminary version of the research \cite{williamsrogers2023designing} by providing new formulations, results, and a more detailed description of the materials.

The remainder of this paper is organized as follows. Section~\ref{sec:lit_rev} reviews state-of-the-art constellation configuration design methods. Section~\ref{sec:mathematical_modeling} introduces the optimization parameters and the five MILP formulations, each illustrated with a corresponding example. Section~\ref{sec:case_studies} presents a comparative analysis between the formulations, exhibits the new formulations derived from the mix and match operations, and the extensions to the proposed formulations. Finally, Sec.~\ref{sec:conclusions} summarizes the main findings and outlines directions for future research.

\section{Literature Review}\label{sec:lit_rev}
This section reviews the literature on satellite constellation configuration design. For each paper, we describe its methodology, the coverage objective (\textit{e.g.,} global, regional, continuous, discontinuous), the fundamental assumptions, and the optimization approach taken, as applicable. Although not comprehensive, this review outlines a wide range of approaches presented in the literature to tackle this complex problem.

A traditional methodology for constellation configuration design involves imposing specific constraints that govern the orbital distribution of satellites. For example, Walker \cite{walker_1970,walker_1977} and Rosette \cite{ballard_taes_1980} patterns enforce all satellites to circular orbits with uniform semi-major axis and inclination, and with a right ascension of the ascending node (RAAN) uniformly distributed. Similarly, Ref.~\cite{luders_aiaa_1961} introduces constellation configurations where the satellites' circular orbits share the same semi-major axis, have polar inclinations, and a uniform RAAN distribution. 
Reference~\cite{avedano_cmda_2013} proposes the 2D Lattice Flower Constellation (LFC) theory, which enables to construct constellation configurations that are time invariant, either with respect to an inertial or rotating reference frames, and whose satellites are uniformly distributed, sharing circular or critically inclined orbits with uniform same semi-major axis and inclination. Reference~\cite{davis_cmda_2013} presents the 3D LFC theory, a generalization of the 2D LFC theory to elliptical orbits of common semi-major axis, inclination, and eccentricity. In general, these highly constrained constellation configurations have a straightforward design approach given the reduced degrees of freedom in their design space, and are proven effective for global coverage with multi-fold time-invariant coverage requirements. However, as a consequence of these constraints, their application to mission scenarios with complex coverage requirements (\textit{e.g.,} multi-fold time-varying) is challenging \cite{lee2020satellite}.

Enumeration procedures, that is, a comprehensive evaluation of the solution space, are employed to design constellation configurations that impose strict orbital distribution requirements. For instance, the authors of Refs.~\cite{walker_1970,ballard_taes_1980} minimize the number of satellites required in the constellation to achieve continuous global coverage. Reference~\cite{ulybyshev_jsr_1999} minimizes the number of satellites required to deliver multi-fold global coverage, adopting streets-of coverage and enforcing all orbits to be circular. Furthermore, enumeration procedures are used to design symmetrical constellation configurations adopting multi-fold regional coverage as the key figure of merit and enforcing circular \cite{adams_jas_1987}, polar \cite{beste_taes_1978}, and elliptical orbits \cite{draim_jgcd_1985}. In addition to adopting regional coverage, Ref.~\cite{lin_jsr_2005} designs a Walker Delta constellation imposing strict MRT mission requirements. The adoption of enumeration procedures provides a simple yet effective design methodology for mission scenarios embracing prescribed constellation configurations; however, their application is generally limited to small mission scenarios given the required full solution space enumeration.

Motivated by the limitations of the enumeration procedures, the literature proposes optimization methods with sophisticated solution space exploration. For example, genetic algorithms (GA) are leveraged to minimize the constellation configuration cost and maximize its robustness against satellite failure, such that continuous global coverage is achieved while enforcing a two-layer Walker Delta pattern \cite{landsard_aa_1998}, and Ref.~\cite{meziane} minimizes the number of satellites and their altitude seeking to achieve continuous regional coverage requiring circular orbits with uniform inclination and RAAN distribution. Further, Ref.~\cite{mortari_ieee_2011} maximizes the global percentage coverage, enforcing either a Walker Delta pattern or a flower constellation. Aiming to gain additional degrees of freedom in the design space, the literature relaxes the requirement of adopting prescribed constellation configurations. As an illustration, Ref.~\cite{wil2001average} concurrently minimizes the ART and MRT over specific regions of interest, allowing both symmetrical and asymmetrical constellation configurations. Variations of the GA algorithm are additionally used to tackle the constellation configuration design problem. For instance, the non-dominated sorting genetic algorithm II is proposed for minimizing, in addition to other objectives, the MRT over specific regions of interest \cite{ferringer2006satellite,ferringer2007efficient}. Conversely to metaheuristics, Refs.~\cite{lee2020satellite,lee2023mclp} tackle the problem of constellation configuration design as an ILP, minimizing the constellation configuration cost to deliver regional complex coverage, and maximizing the observation rewards and percentage coverage, respectively. In addition, alternative optimization methods adopted for global coverage include the gradient descent algorithm \cite{capez_taes_2022} and heuristics \cite{wang_wcl_2022}, and nonlinear models derived from the coverage geometry tackled with commercial solvers to deliver regional coverage \cite{chadalavada}.

Significant efforts have been made in the literature to develop design (optimization) methods for determining optimal satellite constellation configurations. However, to the best of the authors’ knowledge, existing approaches remain highly problem-specific and lack coherence. Most studies focus on a particular coverage-related figure of merit, leading to differing fundamental assumptions, modeling strategies, and solution techniques. This fragmentation highlights a critical gap in the literature and serves as the primary motivation for this work. In this paper, we propose a collection of MILP formulations that we believe are both practically critical and sufficiently distinct to warrant recognition as a unique formulation in their most basic form, as evidenced by the existing literature. The goal of the proposed collection is to provide users with a go-to reference for easy implementation via commercial off-the-shelf solvers, while offering a desirable property (\textit{i.e.}, certificate of optimality).

\section{Constellation Configuration Design Optimization Problem Formulations}\label{sec:mathematical_modeling}
This section introduces the five formulations of the MILP collection, where each one is accompanied by an illustrative example solved to optimality. The formulations are:
\begin{enumerate}
        \item SCLP \cite{lee2020satellite} (Sec.~\ref{sec:sclp}): Determines the minimum-cost constellation configuration that satisfies the complex spatiotemporal coverage requirements.
        \item PSCLP (Sec.~\ref{sec:psclp}): Determines the minimum-cost constellation configuration that satisfies the percentage coverage requirements.
        \item MCLP \cite{lee2023mclp} (Sec.~\ref{sec:mclp}): Determines the constellation configuration, with a user-defined number of satellites, that maximizes the observation rewards.
        \item MMRT (Sec.~\ref{sec:mmrt}): Determines the constellation configuration, with a user-defined number of satellites, that minimizes the MRT.
        \item MART (Sec.~\ref{sec:mart}): Determines the constellation configuration, with a user-defined number of satellites, that minimizes the ART.
\end{enumerate}
In addition, general notation and definitions common to all formulations are provided.

\subsection{Notation and Definitions}\label{sec:not_def}
We define the following general sets, parameters, and decision variables applicable to all formulations discussed in this paper. A summary is provided in Table~\ref{table:parameters_general}. Definitions specific to individual formulations will be addressed in their respective subsections.

\begin{table}[htbp]
\centering
\caption{General sets, parameters, and decision variables for the MILP formulations.}
\begin{tabular}{lll}
\hline
\hline
Type & Symbol & Description \\
\hline
Sets    &$\mathcal{T}$ & Set of time steps (index $t$; cardinality $T$) \\
        &$\mathcal{J}$ & Set of orbital slots (index $j$; cardinality $J$) \\
        & $\mathcal{P}$ & Set of targets (index $p$; cardinality $P$) \\
Parameters     &$c_{j}$ & Cost of orbital slot $j$\\
                &$r_{tp}$ & Coverage threshold for target $p$ at time step $t$\\
                &$b_{tp}$ & Number of satellites visible to target $p$ at time step $t$\\
                &$V_{tjp}$ & $\begin{cases}
                  1, & \text{if orbital slot $j$ is visible from target $p$ at time step $t$} \\
                  0, & \text{otherwise}
                \end{cases}$ \\
Decision variables &  $x_{j}$ & $\begin{cases}
                     1, &\text{if a satellite occupies orbital slot $j$} \\
                    0, &\text{otherwise}
                    \end{cases}$ \\
\hline
\hline
\label{table:parameters_general}
\end{tabular}
\end{table}

Let $\mathcal{T}$ denote the set of discrete time steps, with index $t$ and cardinality $T$. The definition of this set depends on the mission epoch, mission duration, and time step size (or the number of time steps). Let $\mathcal{J}$ denote the set of orbital slots, with index $j$ and cardinality $J$. Each orbital slot is characterized by a unique set of orbital elements (or a state vector) defined at the epoch, without requiring it to conform to a prescribed orbital geometry or satellite distribution rule (\textit{e.g.}, Walker). Associated with each orbital slot is a cost parameter $c_{j} \in \mathbb{R}_{\ge0}$. This cost can be defined differently for various mission scenarios, aiming to capture specific features. For instance, $c_{j}$ could represent the cost of deploying a satellite to orbital slot $j$, or it could encode the station-keeping cost associated with that orbital slot during the mission. In addition, the constellation configuration is defined through the constellation pattern vector $\bm{x}:=\{x_{j} \in \{0,1\}: j \in \mathcal{J}\}$, where each of its elements is given as:
 \begin{equation}
    x_{j} = \begin{cases}
        1, &\text{if a satellite occupies orbital slot $j$}\\
        0, &\text{otherwise}
    \end{cases}
\end{equation}

Let $\mathcal{P}$ denote the set of targets, with index $p$ and cardinality $P$. At time step $t$, target $p$ has associated spatial coordinates (\textit{e.g.,} geodetic, Cartesian), which enable the definition of static or dynamic targets (\textit{e.g.,} ground-, space-based), visibility requirements (\textit{e.g.,} minimum-elevation angle), and a coverage threshold parameter $r_{tp} \in \mathbb{Z}_{> 0}$. In particular, we define target $p$ as covered if it is visible to at least $r_{tp}$ satellites at time step $t$, allowing users to specify multi-fold, time-varying coverage requirements throughout the entire mission or within specific time windows.

The visibility state of orbital slot $j$ at time step $t$ over target $p$ is encoded using Boolean visibility state parameter $V_{tjp} \in \{0,1\}$, defined as:
\begin{equation}
    V_{tjp} = \begin{cases}
        1, &\text{if target $p$ is visible from orbital slot $j$ at time step $t$}\\
        0,  &\text{otherwise}
    \end{cases}
\end{equation}
To construct $V_{tjp}$, we propagate orbital slot $j$ throughout the mission horizon, we apply visibility masking and encode the corresponding Boolean visibility state $V_{tjp}$ over target $p$ for each time step $t$. The orbit propagation and the visibility computation can be implemented by adopting a custom or commercial off-the-shelf propagator (\textit{e.g.,} STK's \texttt{HPOP} \cite{STK}, MATLAB's \texttt{propagateOrbit} \cite{MATLAB}) and satellite-target access computation algorithms (\textit{e.g.,} MATLAB's \texttt{access} \cite{MATLAB}), respectively.

We define the visibility state over target $p$ for a given constellation pattern $\bm{x}$ leveraging the constellation coverage timeline $\bm{b}_{p}:=\{b_{tp} \in \mathbb{Z}_{\ge0}: t \in \mathcal{T}\}$, where each of its elements indicate the number of satellites visible to target $p$ at time step $t$, and is given as \cite{lee2023mclp}:
\begin{equation}
    b_{tp} = \sum_{j \in \mathcal{J}}V_{tjp}x_{j}\label{eq:cov_timeline_opt}
\end{equation}

\begin{figure}[htpb]
    \centering
    \includegraphics[width=0.6\linewidth]{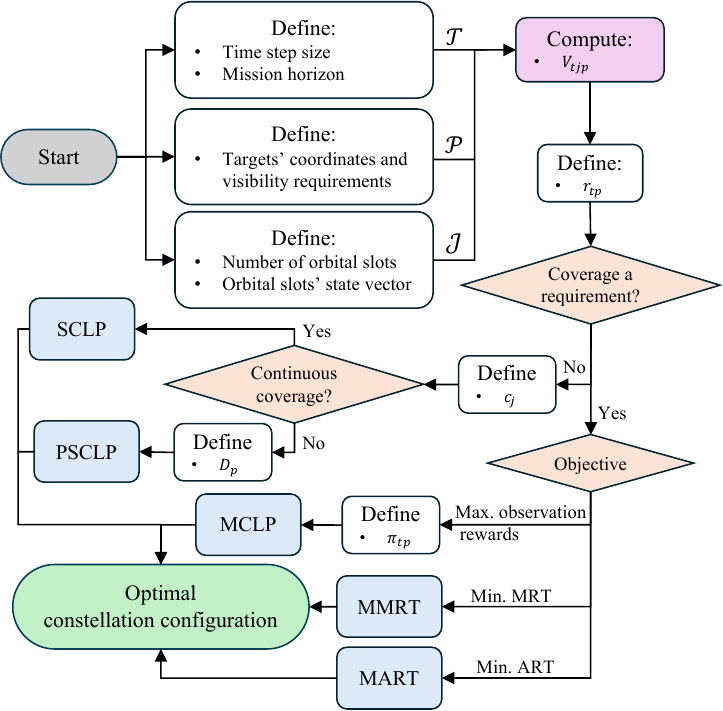}
    \caption{Flowchart for the optimal constellation configuration design.}
    \label{fig:flowchart}
\end{figure}

Figure~\ref{fig:flowchart} presents the proposed constellation configuration design optimization flowchart. In particular, the definition of $\mathcal{T}$, $\mathcal{P}$, and $\mathcal{J}$ stems from a trade-off between the mission requirements (\textit{e.g.,} mission horizon, targets of interest, available orbits to locate satellites), the fidelity of the model (\textit{e.g.,} required time step size to accurately compute coverage considering sensor specifications), and the size of the problem, which impacts the scalability of the optimization.

\subsection{Common Parameters in Illustrative Examples} \label{sec:parameters}
We present each formulation accompanied by an illustrative example. All illustrative examples share the same set of parameters and are solved to optimality using the Gurobi Optimizer version 12.0.2 with default settings on an Intel Core i9-14900 at 2.00 GHz base speed. Without loss of generality, the parameters are selected for ease of exposition of each formulation's unique features and to facilitate intuitive graphical and analytical contrast between the formulations.

Each illustrative example has 716 time steps of size \SI{120}{s}, with an epoch defined as January 1, 2025, at 12:00:00.000 Universal Coordinated Time (UTC). Seattle, Washington, USA (geodetic coordinates \{$47.36\degree$N, $122.19\degree$W, \SI{0}{m}\}) is designated as the sole target with the coverage threshold parameter set as $r_{tp}=1$ for all $t \in \mathcal{T}$, that is, continuous single-fold coverage, and a minimum elevation angle of \SI{5}{deg}. 

All orbital slots are placed in circular sun-synchronous common repeating ground track (RGT) orbits with a semi-major axis of \SI{11546.32}{km}, an inclination of \SI{142.14}{deg}, and a 7:1 resonance ratio (\textit{i.e.}, satellites complete seven ascending node crossings during one nodal period of Greenwich \cite{lee2020satellite,mortari_jas_2004,viptil_jsr_2012}), yielding a repetition period of \SI{85953.92}{s}. The adoption of common RGT orbits enables us to present the results as prescribed by the access-pattern-coverage (APC) decomposition~\cite{lee2020satellite}. The APC decomposition is a constellation-coverage model that leverages three elements. The first element is target $p$'s reference visibility profile $\bm{v}_{p}:=\{v_{tp}\in \{0,1\}:t\in \mathcal{T}\}$, where each $v_{tp}$ encodes the reference orbital slot's visibility state at time step $t$. The second element is the constellation pattern vector $\bm{x}$, where each of its elements $x_j$ is obtained by temporally shifting the reference orbital slot along the common ground track. The third element is target $p$'s coverage timeline $\bm{b}_{p}$, obtained as a circular convolution between the reference visibility profile $\bm{v}_{p}$ and the constellation pattern vector $\bm{x}$. We refer the reader to Refs.~\cite{lee2020satellite,lee2023mclp} for a comprehensive definition of the APC decomposition, and to Appendix~\ref{app:oe_generation} for a definition of RGT orbital slots generation.

\subsection{Set Covering Location Problem}\label{sec:sclp}
SCLP determines the minimum-cost constellation configuration that strictly satisfies the complex spatiotemporal coverage requirements. This formulation, originally proposed in Ref.~\cite{lee2020satellite} for satellite constellation configuration design, resembles the SCLP formulation \cite{toregas_or_1971}, a class of FLP.

SCLP minimizes the cost of the constellation configuration encoded in objective function:
\begin{equation}
    \sum_{j\in \mathcal{J}}c_{j} x_{j} \label{sclp:obj}
\end{equation}

For each time step $t$, we enforce the spatiotemporal coverage requirement, that is, at least $r_{tp}$ satellites visible to target $p$ with constraints:
\begin{equation}
     \sum_{j\in \mathcal{J}} V_{tjp}x_{j} \ge r_{tp}, \quad  \forall t \in \mathcal{T}, \ \forall p \in \mathcal{P} \label{sclp:visibility}
\end{equation}

Piecing it all together, the SCLP formulation is as follows \cite{lee2020satellite}:
\begin{alignat}{2}
\textsc{(\hypertarget{sclp}{SCLP})} \quad \min \quad & \sum_{j\in \mathcal{J}}c_{j} x_{j} \notag \\
\text{s.t.} \quad & \sum_{j\in \mathcal{J}} V_{tjp}x_{j} \ge r_{tp}, &\quad  \forall t \in \mathcal{T}, \ \forall p \in \mathcal{P} \notag \\
&x_{j} \in \{0, 1\}, &\quad \forall j\in \mathcal{J} \label{sclp:xj}
\end{alignat} 

\begin{figure}[htbp]
	\centering
	\begin{subfigure}[c]{0.492\linewidth}
		\centering
		\includegraphics[width=\linewidth]{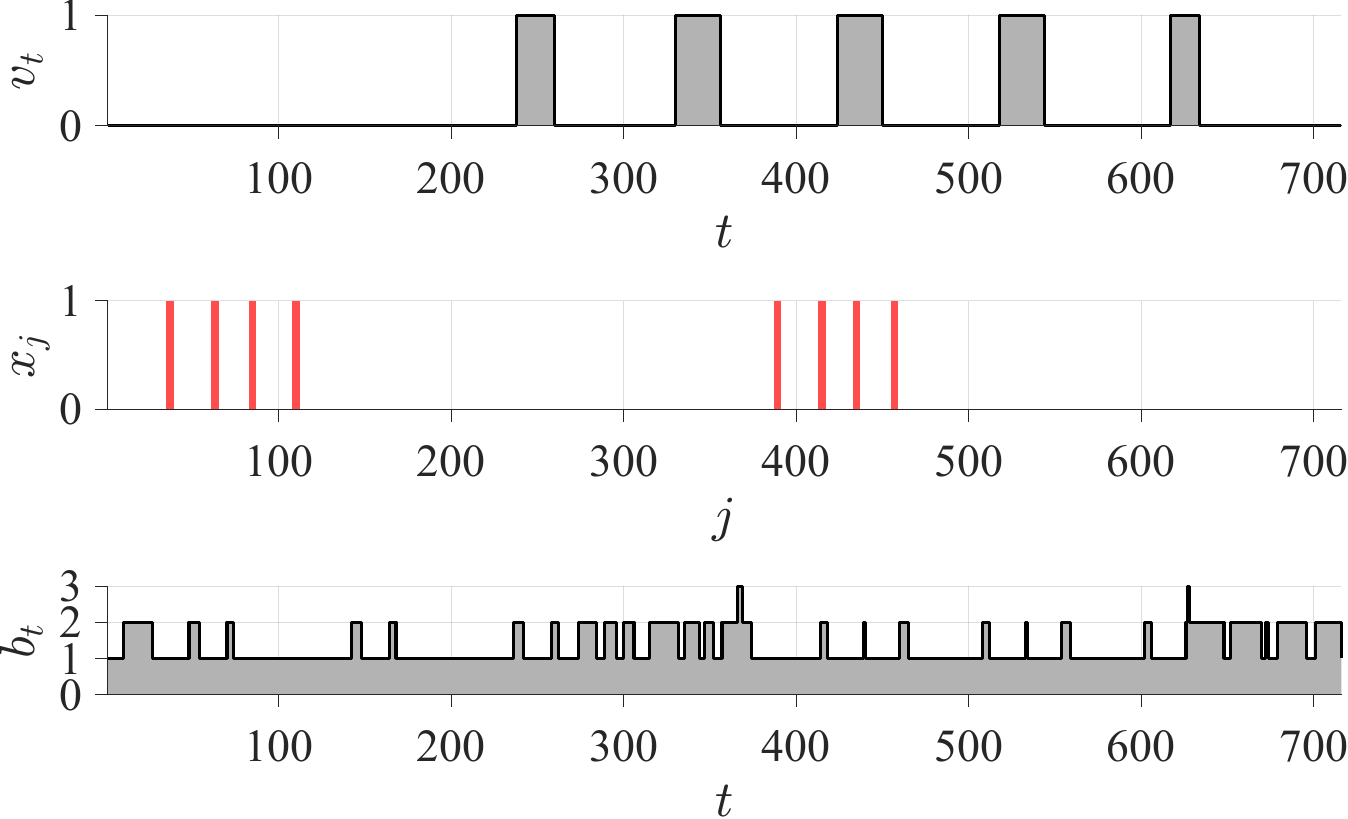}
		\caption{Illustration of the APC decomposition.}
            \label{fig:sclp_a}
	\end{subfigure}
	\begin{subfigure}[c]{0.492\linewidth}
		\centering
		\includegraphics[width=\linewidth]{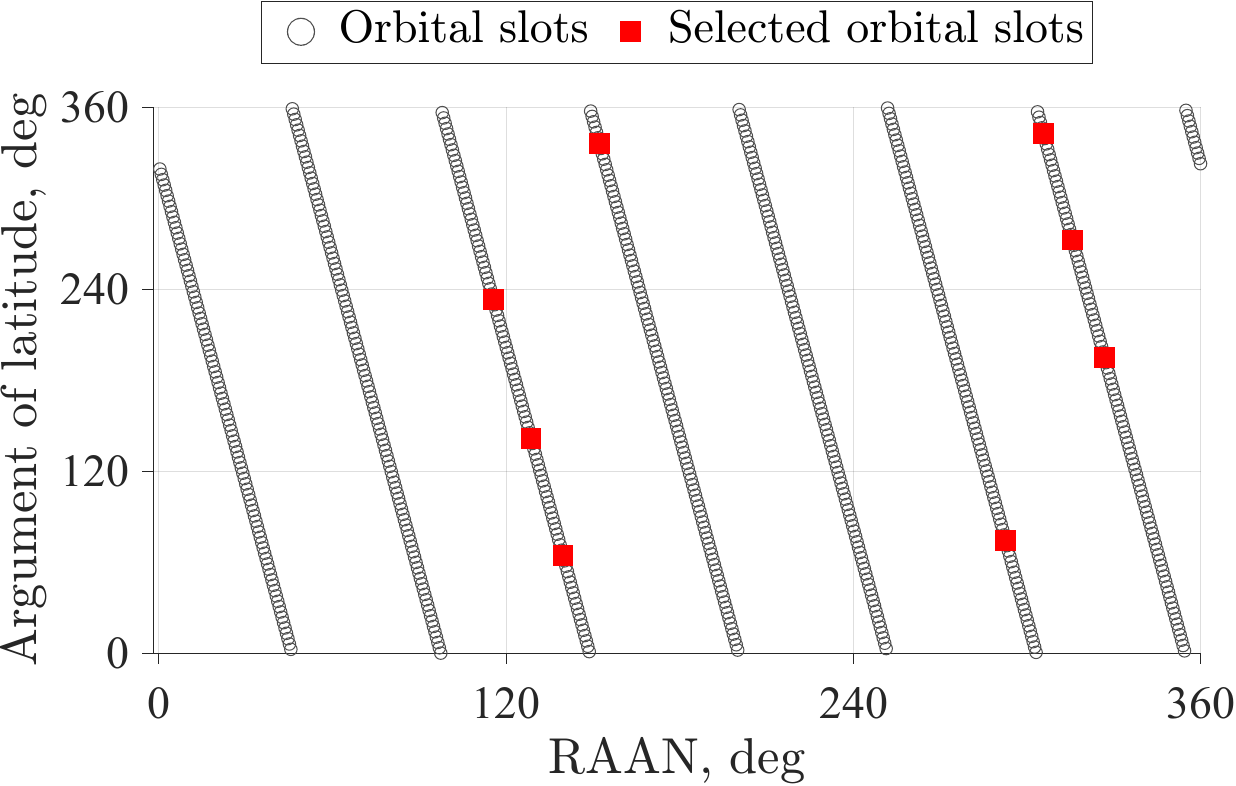}
		\caption{Optimal constellation configuration's orbital slot distribution.}
            \label{fig:sclp_b}
	\end{subfigure}
        \begin{subfigure}[c]{0.5853\linewidth}
            \centering
             \includegraphics[width=0.60\linewidth]{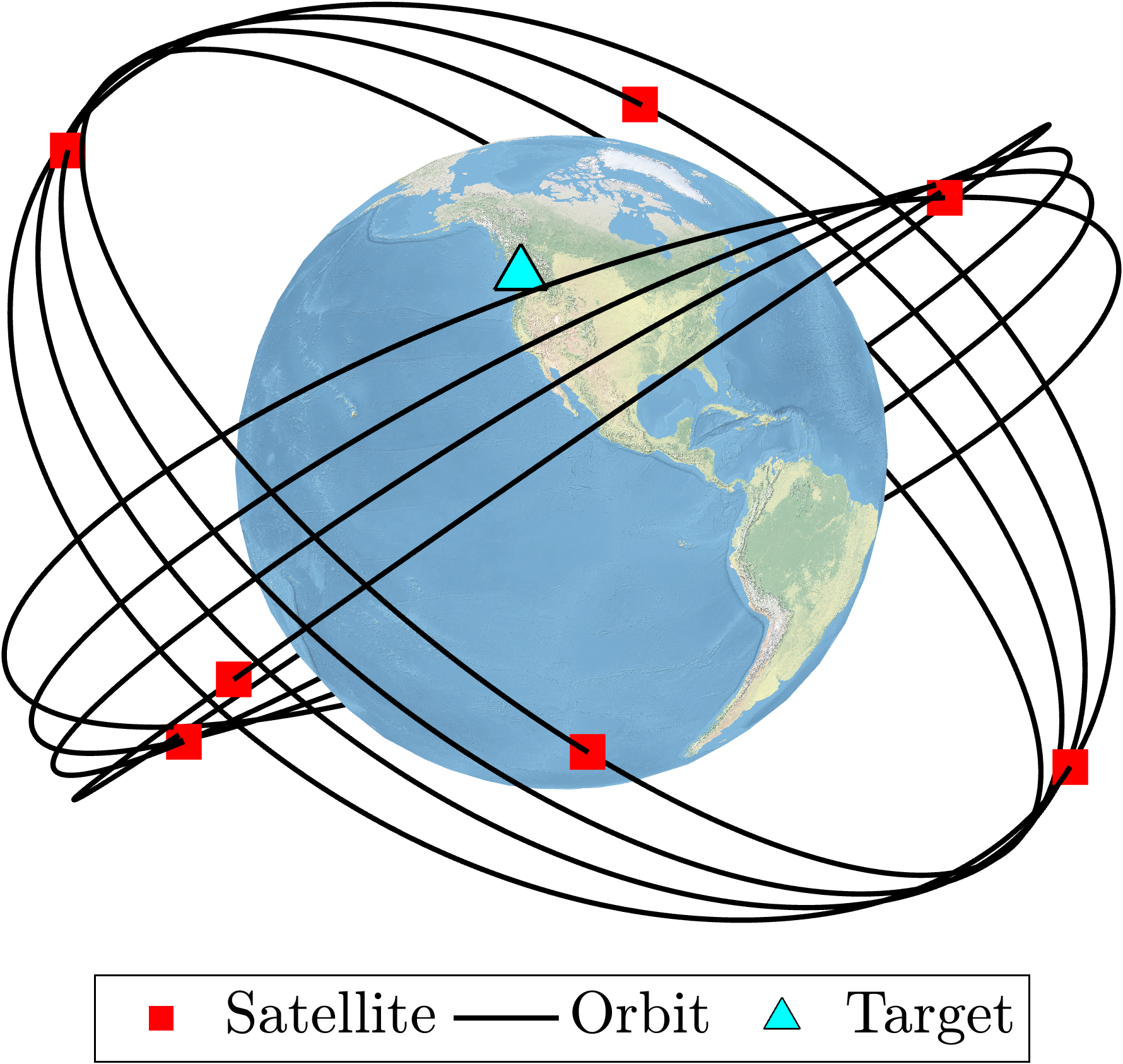}
            \caption{Eight-satellite optimal constellation configuration.}
                \label{fig:sclp_c}
        \end{subfigure}
	\caption{SCLP formulation results.}
	\label{fig:sclp}
\end{figure}

\begin{example}[Single-Fold Continuous Coverage SCLP] \label{eg:sclp}
Given a cost $c_j=1$ for all $j \in \mathcal{J}$, the optimal objective value is eight, and the runtime is \SI{16.13}{s}, corresponding to an optimal constellation configuration of eight satellites. Figure~\ref{fig:sclp_a} showcases the APC decomposition \cite{lee2020satellite} of the optimal constellation. At the top of it lies the reference visibility profile $\bm{v}_{p}$, at the middle the constellation pattern vector $\bm{x}$, indicating the occupied orbital slots with red impulses, and at the bottom the coverage timeline $\bm{b}_{p}$. Figure~\ref{fig:sclp_b} illustrates the distribution of all orbital slots in the RAAN versus argument of latitude plane, with the occupied orbital slots indicated with red squares. Lastly, Fig.~\ref{fig:sclp_c} presents a 3D visualization of the optimal constellation configuration with the satellites' orbits in the Earth-centered inertial (ECI) frame, and the target at the epoch.
\end{example}

\subsection{Partial Set Covering Location Problem}\label{sec:psclp}
PSCLP determines the minimum-cost constellation configuration that satisfies each target's minimum temporal percentage coverage requirement. The PSCLP formulation holds practical importance as it allows mission designers to implement the concept of temporal percentage coverage, a universally adopted figure of merit in constellation design and analysis \cite{SMAD}. This formulation, similarly to SCLP, resembles the PSCLP class of FLP \cite{daskin_ga_1999}. Further, PSCLP generalizes SCLP by incorporating more flexible coverage requirements. However, unlike SCLP, which seeks to fully satisfy the complex coverage requirements imposed on all targets, PSCLP permits partial fulfillment, hence percentage coverage. Table~\ref{table:parameters_psclp} presents additional parameters and variables used in PSCLP.

\begin{table}[htbp]
\caption{Specific PSCLP parameters and decision variables.}
\centering
\begin{tabular}{lll}
\hline
\hline
Type & Symbol & Description \\
\hline
Parameters  &$D_{p}$ & Temporal percentage coverage for target $p$\\
            &$D$ & Mean spatiotemporal percentage coverage\\
Decision variables &  $y_{tp}$ & $\begin{cases}
                     1, &\text{if target $p$ is covered at time step $t$} \\
                    0, &\text{otherwise}
                    \end{cases}$ \\
\hline
\hline
\label{table:parameters_psclp}
\end{tabular}
\end{table}

PSCLP minimizes the total cost of the constellation configuration leveraging objective function~\eqref{sclp:obj}. Further, to extend the temporal percentage coverage requirements, PSCLP uses the coverage state decision variables $\bm{y}:=\{y_{tp} \in \{0,1\}:t \in \mathcal{T}, p \in \mathcal{P}\}$, where each element $y_{tp}$ is defined as:
\begin{equation}
  y_{tp} = \begin{cases}
    1, & \text{if target $p$ is covered at time step $t$} \\
    0, & \text{otherwise}
\end{cases}  
\end{equation}

For each time step $t$, target $p$'s coverage state decision variable $y_{tp}$ is activated (\textit{i.e.,} equal to one), if at least $r_{tp}$ satellites are visible to it. Then, to enforce this condition, we propose constraints:
\begin{equation}
    \sum_{j\in\mathcal{J}} V_{tjp}x_{j}\ge r_{tp} y_{tp}, \quad \forall t\in \mathcal{T}, \forall p \in\mathcal{P} \label{psclp:V_x_y}\\
\end{equation}
Then, we enforce target $p$'s minimum temporal percentage coverage requirement $D_{p} \in \mathbb{Z}_{\ge0}$, that is, the sum of the covered time steps for target $p$, with constraints:
\begin{equation}
    \sum_{t\in \mathcal{T}} y_{tp} \ge D_{p},  \quad \forall p \in\mathcal{P} \label{psclp:y_dp} 
\end{equation}

Lastly, the PSCLP formulation is given as:
\begin{alignat}{2}
    \textsc{(\hypertarget{psclp}{PSCLP})} \quad \min \quad & \sum_{j\in \mathcal{J}}c_{j} x_{j} \notag \\
    \text{s.t.} \quad & \sum_{j\in\mathcal{J}} V_{tjp}x_{j}\ge r_{tp} y_{tp}, &\quad \forall t\in \mathcal{T}, \forall p \in\mathcal{P} \notag\\
    & \sum_{t\in \mathcal{T}} y_{tp} \ge D_{p}, & \quad \forall p \in\mathcal{P} \notag\\
    &   x_{j} \in \{0, 1\}, &\quad \forall j\in \mathcal{J} \notag\\
    & y_{tp} \in \{0,1\}, &\quad  \forall t\in \mathcal{T}, \forall p \in\mathcal{P} \label{psclp:ytp} 
\end{alignat}

\begin{remark}[PSCLP as a Generalization of SCLP]
    It should be noted that if $D_{p} = T$, the temporal percentage coverage requirement of target $p$ corresponds to continuous coverage. Therefore, PSCLP determines the minimum-cost constellation configuration that delivers continuous coverage over each target $p$. Hence, it resembles SCLP.
\end{remark}

\begin{remark}[Mean Percent Coverage]
    A variation of this formulation is possible if, instead of setting individual temporal percent coverage requirements as in constraints~\eqref{psclp:y_dp} per target $p$, one may impose a \textit
    {mean spatiotemporal percent coverage} requirement $D  \in \mathbb{Z}_{\ge 0}$. In this case, constraints~\eqref{psclp:y_dp} are replaced by the following constraint:
\begin{equation}
\sum_{p\in\mathcal{P}}\sum_{t\in\mathcal{T}}y_{tp}\ge D
\end{equation}
\end{remark}

\begin{figure}[htbp]
	\centering
	\begin{subfigure}[c]{0.492\linewidth}
		\centering
		\includegraphics[width=\linewidth]{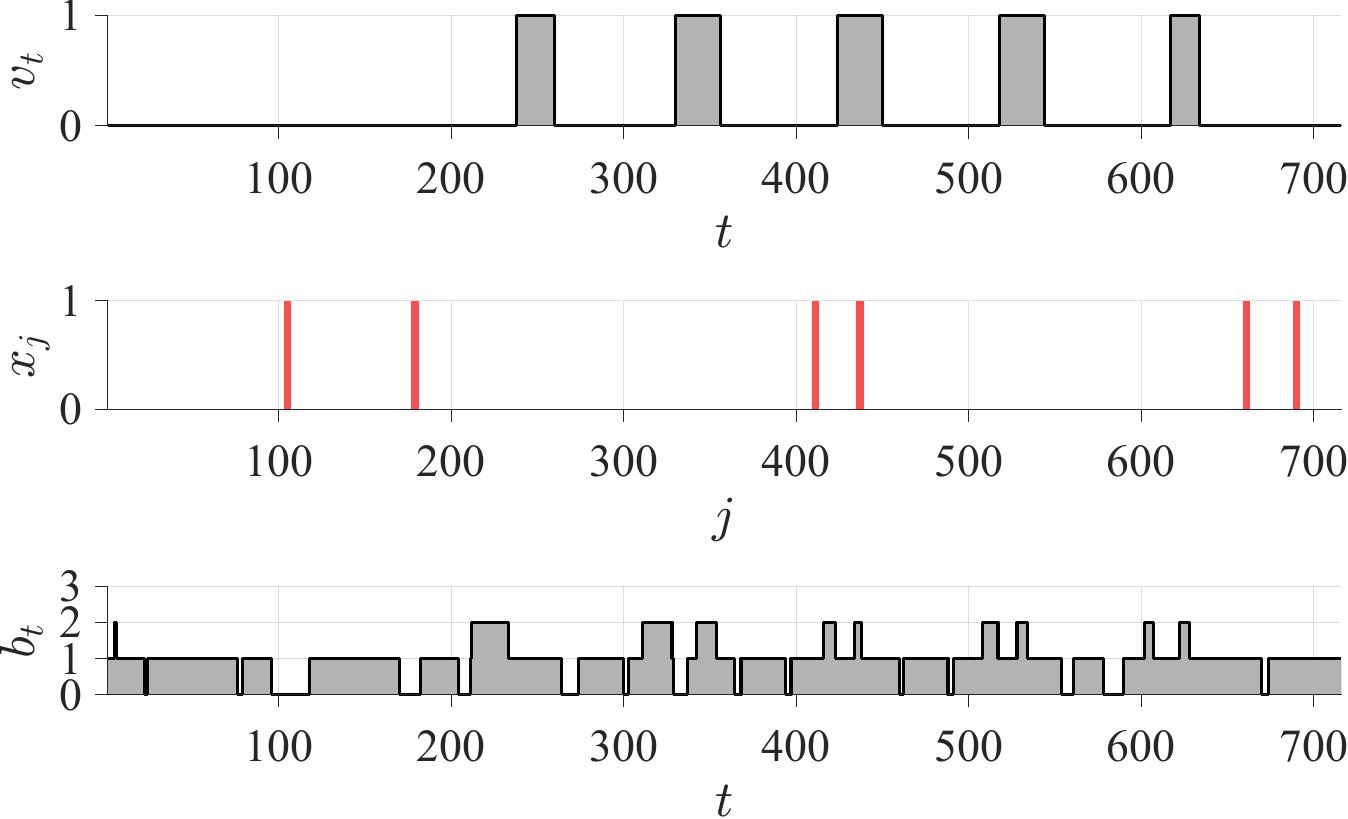}
		\caption{Illustration of the APC decomposition.}
            \label{fig:psclp_a}
	\end{subfigure}
	\begin{subfigure}[c]{0.492\linewidth}
		\centering
		\includegraphics[width=\linewidth]{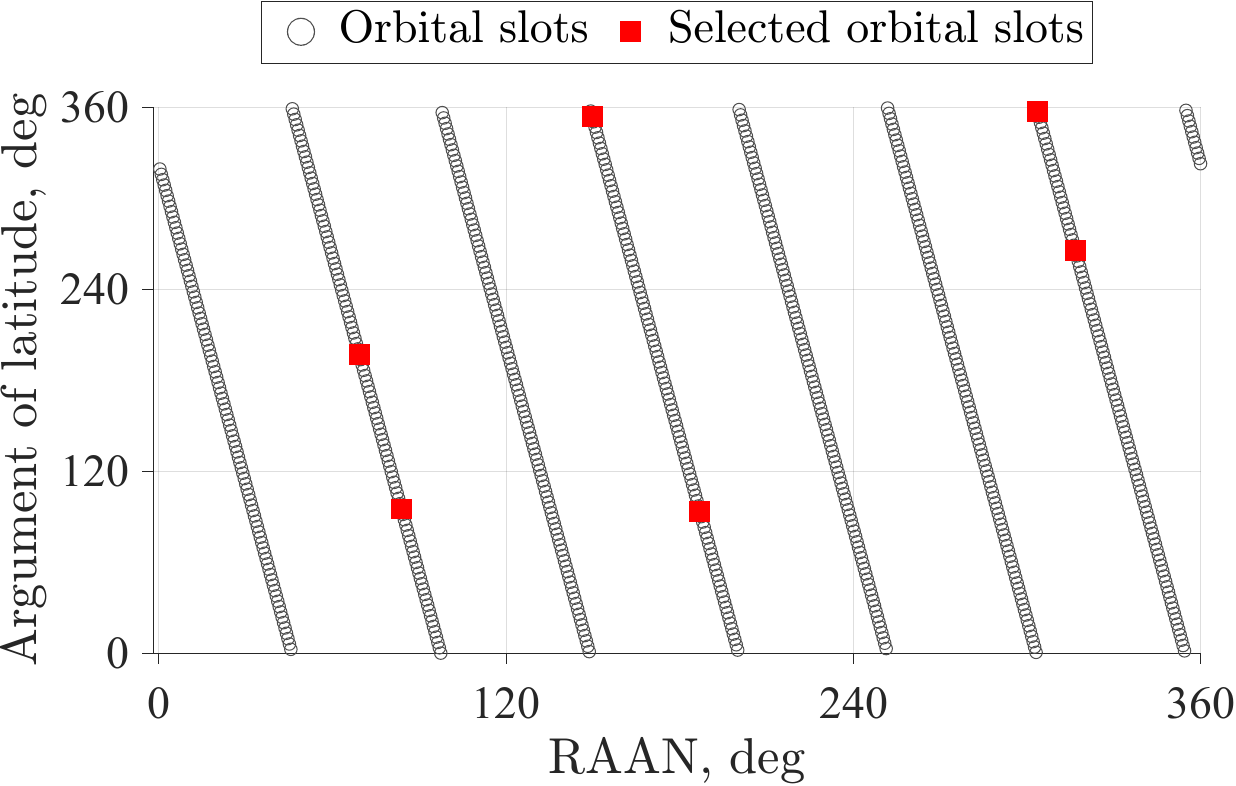}
		\caption{Optimal constellation configuration's orbital slot distribution.}
            \label{fig:psclp_b}
	\end{subfigure}
    \begin{subfigure}[c]{0.5853\linewidth}
            \centering
             \includegraphics[width=0.60\linewidth]{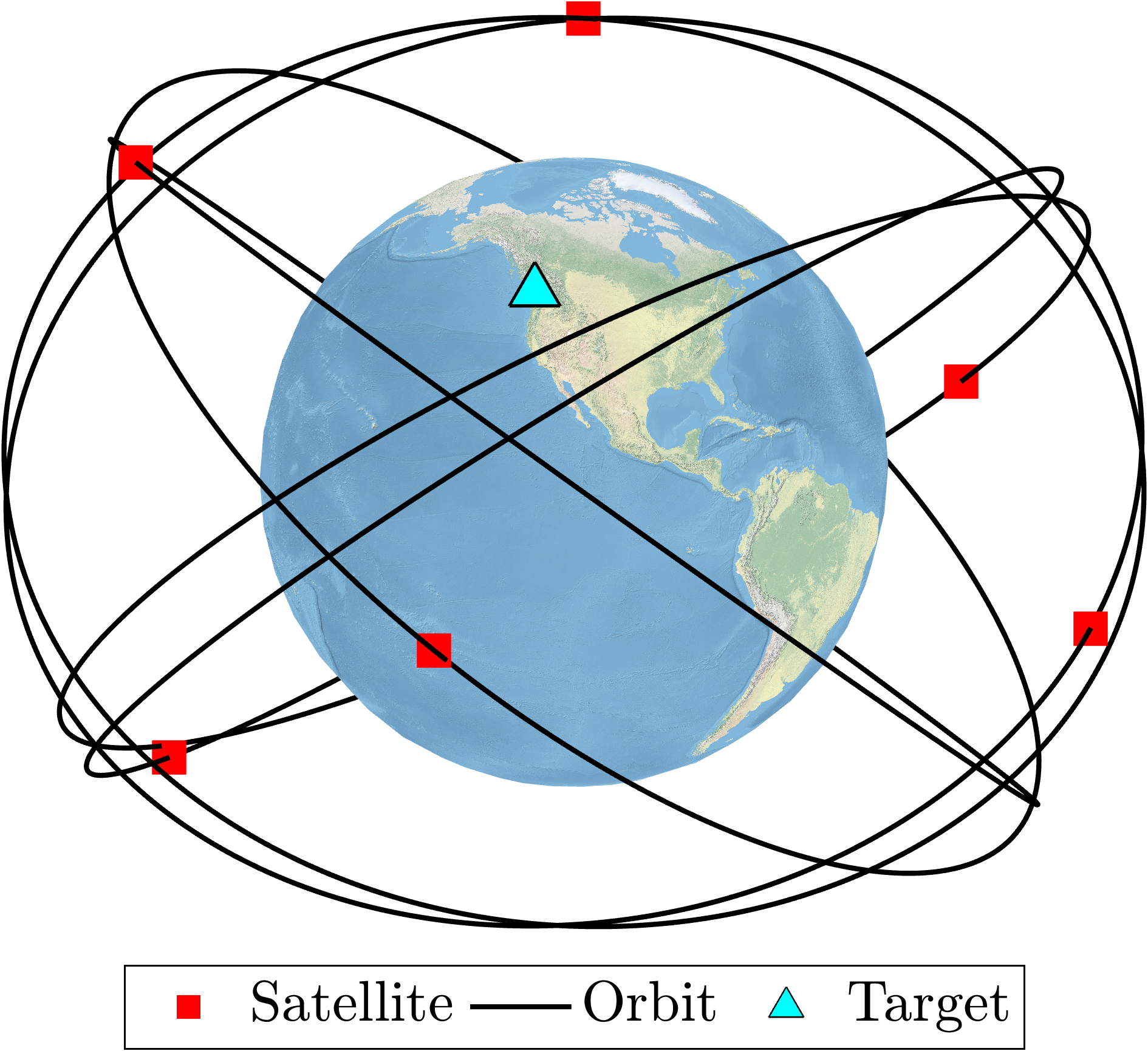}
            \caption{Six-satellite optimal constellation configuration with 85.75 \% coverage.}
                \label{fig:psclp_c}
        \end{subfigure}
	\caption{PSCLP formulation results.}
	\label{fig:psclp}
\end{figure}

\begin{example}[Single-Fold \SI{85}{\%} Coverage PSCLP]
   Given a cost $c_j=1$ for all $j \in \mathcal{J}$, and a required minimum percentage coverage of \SI{85}{\%}, the optimal objective value is six, and the runtime is \SI{3.82}{s}, corresponding to an optimal constellation configuration of six satellites. The optimal constellation configuration delivers \SI{85.75}{\%} coverage over Seattle. Figure~\ref{fig:psclp_a} presents the APC decomposition, and Fig.~\ref{fig:psclp_b} outlines the distribution of the orbital slots in the RAAN versus argument of latitude plane. Figure~\ref{fig:psclp_c} illustrates a 3D visualization of the constellation configuration and the target at the epoch.
\end{example}

\subsection{Maximal Covering Location Problem}\label{sec:mclp}
MCLP determines the location of $N$ satellites within a set $\mathcal{J}$ to maximize the sum of time-dependent observation rewards associated with a target. The significance of MCLP applied to the constellation configuration design optimization, as originally proposed by Ref.~\cite{lee2023mclp}, is that it enables the user to design a constellation where the number of satellites is a parameter and the objective is to maximize the targets' observational rewards. It should be noted that the MCLP was originally proposed in the context of operations research by Ref.~\cite{church1974maximal}. Further, the objectives of SCLP, PSCLP, and MCLP are distinct; the first two aim to minimize the constellation configuration cost, whereas MCLP maximizes the observation rewards given a fixed number of satellites. Lastly, Table~\ref{table:parameters_mclp} presents the newly introduced parameters for the MCLP formulation.

\begin{table}[htbp]
\caption{Specific MCLP parameters.}
\centering
\begin{tabular}{lll}
\hline
\hline
Type & Symbol & Description \\
\hline
Parameters  & $\pi_{tp}$ & Observation reward for target $p$ at time step $t$\\
                & $N$ & Number of satellites in the constellation\\
\hline
\hline
\label{table:parameters_mclp}
\end{tabular}
\end{table}

MCLP maximizes the sum of the collected time-dependent observation rewards over the set of targets $\mathcal{P}$ using objective function:
\begin{equation}
    \sum_{p\in\mathcal{P}}\sum_{t\in \mathcal{T}}\pi_{tp} y_{tp} \label{mclp:obj} \\
\end{equation}
where $\pi_{tp} \in \mathbb{R}_{\ge0}$ denotes the observation reward associated to target $p$ at time step $t$, and is collected if the target is covered. Further, the coverage state of each target $p$ is enforced by constraints~\eqref{psclp:V_x_y}.

The constellation must use $N$ satellites, as imposed by constraint:
\begin{equation}
    \sum_{j\in \mathcal{J}} x_{j} = N \label{mclp:cardinality} \\
\end{equation}

Given the objective function and constraints, the MCLP formulation is:
\begin{alignat*}{2}
    \textsc{(\hypertarget{mclp}{MCLP})} \quad \max \quad & \sum_{p\in\mathcal{P}}\sum_{t\in \mathcal{T}}\pi_{tp} y_{tp} \\
    \text{s.t.} \quad & \sum_{j\in\mathcal{J}} V_{tjp}x_{j}\ge r_{tp} y_{tp}, &\quad \forall t\in \mathcal{T}, \forall p \in\mathcal{P}\\
    &\sum_{j\in \mathcal{J}} x_{j} = N  \\
    &   x_{j} \in \{0, 1\}, &\quad \forall j\in \mathcal{J} \\
    & y_{tp} \in \{0,1\}, &\quad  \forall t\in \mathcal{T}, \forall p \in\mathcal{P} 
\end{alignat*}

\begin{remark}[MCLP with Budget Constraint]
In certain mission design scenarios, it is beneficial to limit the constellation configuration's budget instead of fixing the number of satellites. In such cases, constraint~\eqref{mclp:cardinality} is replaced by constraint:
    \begin{equation}
    \sum_{j \in \mathcal{J}}c_{j} x_{j} \le C
\end{equation}
where $C \in \mathbb{R}_{>0}$ is the constellation configuration's budget.
\end{remark}

\begin{figure}[htbp]
	\centering
	\begin{subfigure}[c]{0.492\linewidth}
		\centering
		\includegraphics[width=\linewidth]{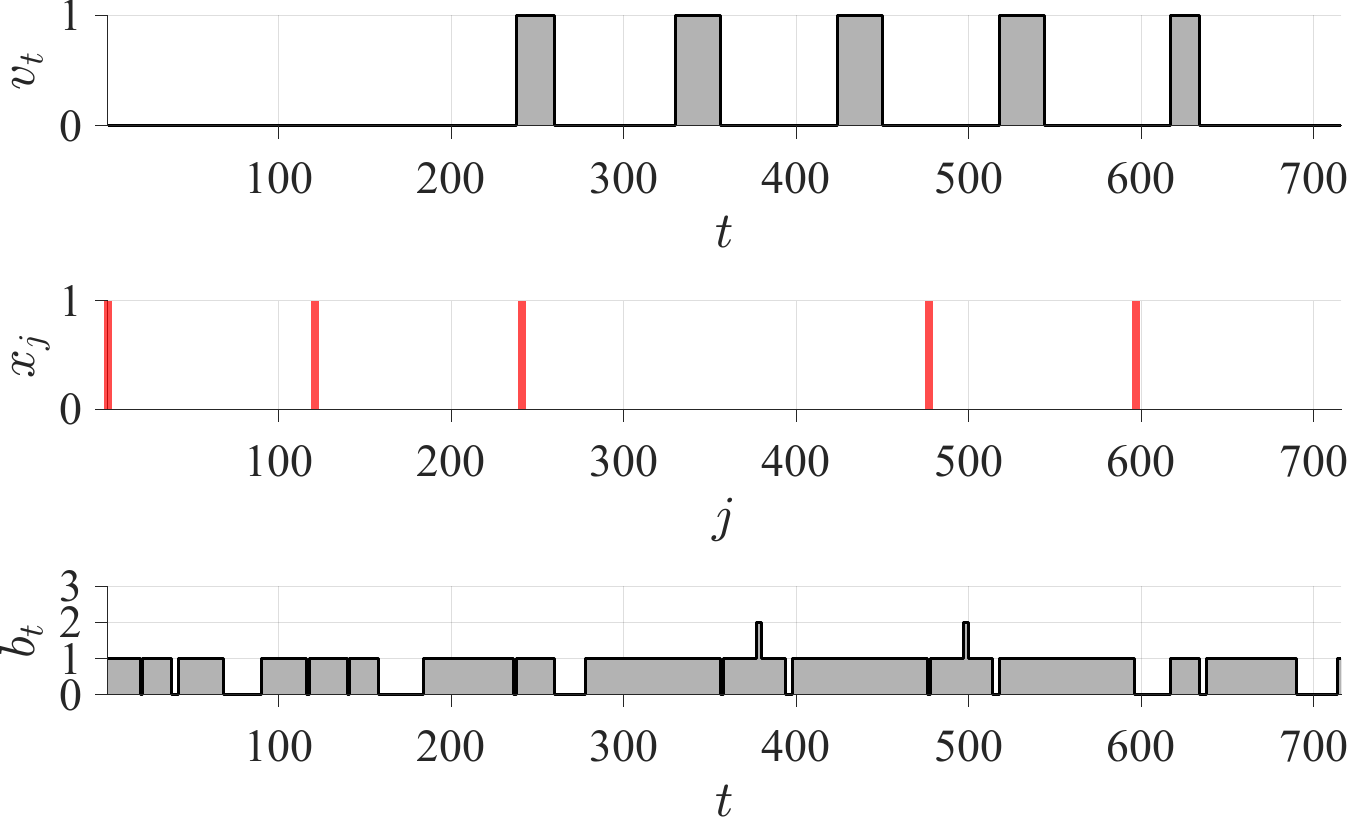}
		\caption{Illustration of the APC decomposition.}
            \label{fig:mclp_a}
	\end{subfigure}
	\begin{subfigure}[c]{0.492\linewidth}
		\centering
		\includegraphics[width=\linewidth]{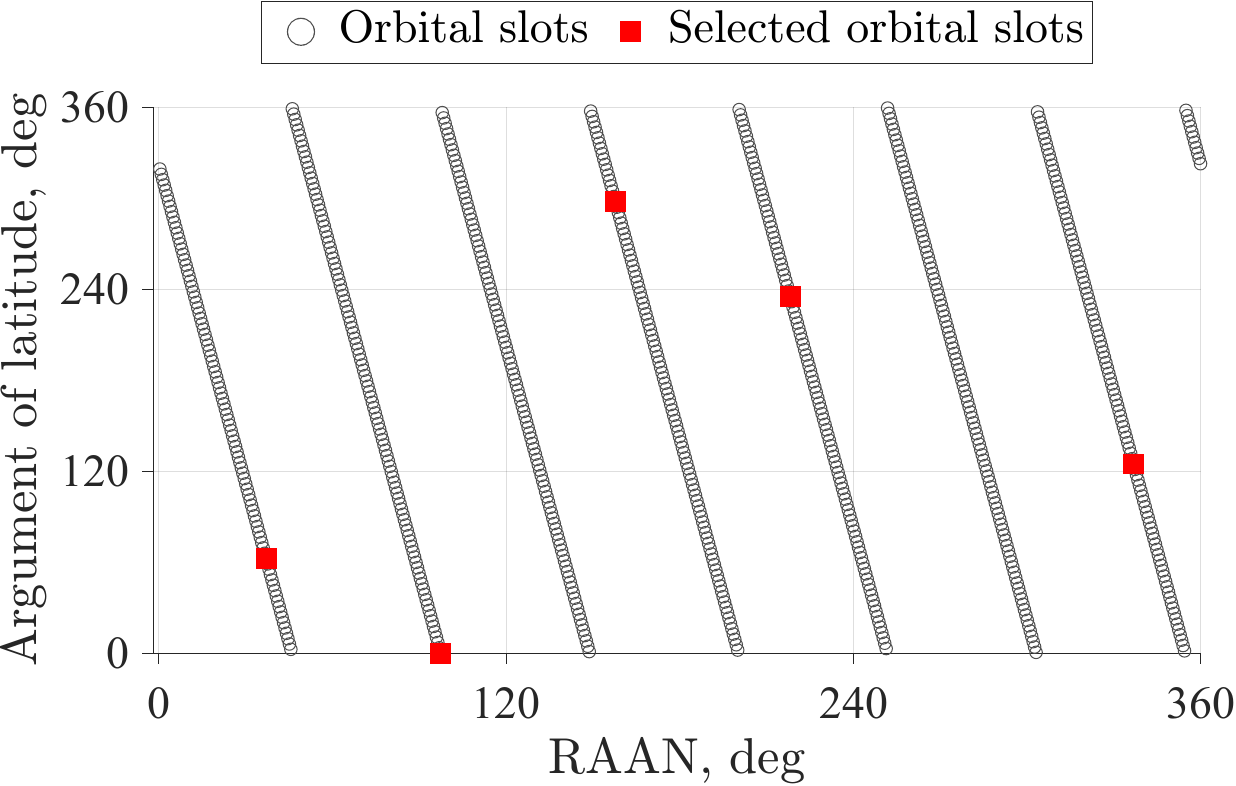}
		\caption{Optimal constellation configuration's orbital slot distribution.}
            \label{fig:mclp_b}
	\end{subfigure}
        \begin{subfigure}[c]{0.5853\linewidth}
            \centering
             \includegraphics[width=0.60\linewidth]{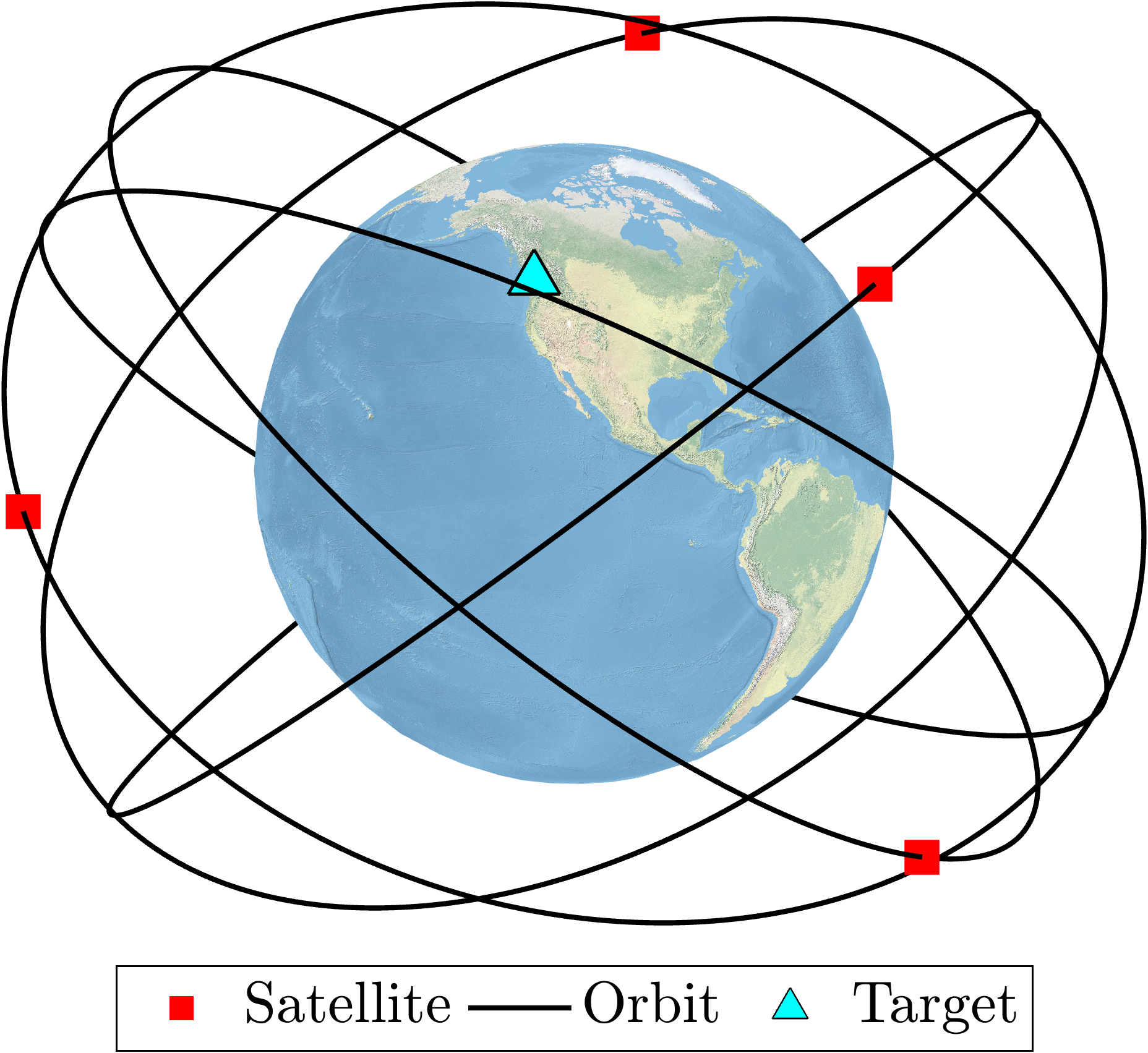}
                \caption{Five-satellite optimal constellation configuration with 80.86 \% coverage.}
                    \label{fig:mclp_c}
            \end{subfigure}
	\caption{MCLP formulation results.}
	\label{fig:mclp}
\end{figure}

\begin{example}[Five-Satellite MCLP]
Given a fixed number of satellites equal to five and an observation reward $\pi_{tp}=1$ for all $t \in \mathcal{T}$, the optimal objective value is 579, and the runtime is \SI{7.88}{s}, corresponding to an optimal constellation configuration that delivers \SI{80.86}{\%} coverage over Seattle. Figure~\ref{fig:mclp_a} presents the APC decomposition, and Fig.~\ref{fig:mclp_b} outlines the distribution of the orbital slots in the RAAN versus argument of latitude plane. Figure~\ref{fig:mclp_c} illustrates a 3D visualization of the constellation configuration and the target at the epoch.
\end{example}

\subsection{Minimal Maximum Revisit Time Problem}\label{sec:mmrt}
MMRT determines the location of $N$ satellites within a set $\mathcal{J}$ to minimize the MRT. The significance of this formulation lies in its ability to enable mission designers to tackle the MRT, a relevant metric for Earth observation missions and telecommunication services, where the maximum length of the coverage gaps impacts the performance of the mission. Although this formulation fixes the number of satellites, it is distinct from MCLP in the sense that it focuses on the MRT instead of maximizing the observation rewards. Table~\ref{table:parameters_mrt} introduces the novel parameter and decision variables used in the formulation.

\begin{table}[htbp]
\caption{Specific MMRT parameters and decision variables.}
\centering
\begin{tabular}{lll}
\hline
\hline
Type & Symbol & Description \\
\hline
Parameter   & $\beta$ & Small positive real-valued auxiliary parameter\\
Decision variables  &  $w_{tp}$ & Duration of coverage gap at time step $t$ for target $p$ \\
      & $Z$ & Nonnegative integer variable encoding the MRT\\
\hline
\hline
\label{table:parameters_mrt}
\end{tabular}
\end{table}

The MMRT minimizes the MRT encoded in decision variable $Z \in \mathbb{Z}_{\ge0}$ using an $N$-satellite constellation configuration, as enforced by cardinality constraints~\eqref{mclp:cardinality}. Conversely to PSCLP and MCLP, MMRT enforces the coverage state for each time step $t$ over each target $p$ using the Big-M method. Then, the visibility constraints are given as:
\begin{subequations}
    \begin{alignat}{2}
    & \sum_{j \in \mathcal{J}}V_{tjp}x_{j} - r_{tp}\ge -J(1-y_{tp}) - \beta, \quad & \forall t \in \mathcal{T}, \forall p \in \mathcal{P}\label{mrt:bigM_1}\\
    & \sum_{j \in \mathcal{J}}V_{tjp}x_{j} - r_{tp}\le Jy_{tp} - \beta, \quad & \forall t \in \mathcal{T}, \forall p \in \mathcal{P}\label{mrt:bigM_2}         
\end{alignat}
\end{subequations}
where $\beta$ is a nonnegative real valued constant lower than one. The motivation for adopting visibility constraints that leverage the Big-M method over constraints~\eqref{sclp:visibility} lies in the different structure of the MMRT formulation. For the MCLP, the optimization aims to maximize the observation rewards by activating coverage state decision variables $y_{tp}$, without imposing any lower bound on their value. Conversely, MMRT aims to minimize the duration of the longest coverage gap; therefore, for specific cases where the target is covered, MMRT is not forced to activate coverage state decision variables $y_{tp}$ equal to one if they do not belong to the longest coverage gap, hence demanding constraints on its lower bound.    
    
We introduce decision variables $w_{tp} \in \mathbb{Z}_{\ge0}$ to count for the number of contiguous time steps where target $p$ is not covered until time step $t$, that is, the duration of the coverage gap. Constraints~\eqref{mrt:w0} initialize decision variables $w_{1,p}$ considering the coverage state at the first time step for each target $p$. If target $p$ is covered at time step $t$, constraints~\eqref{mrt:counter_1} enforce $w_{tp}=0$. Alternatively, if target $p$ is not covered, that is, $y_{tp}=0$, constraints~\eqref{mrt:counter_1} are nonbinding and constraints~\eqref{mrt:counter_2} and~\eqref{mrt:counter_3} yield an equality increasing the value of $w_{tp}$ by one. Further, constraints~\eqref{mrt:z_w} couple MRT decision variable $Z$ with the largest coverage gap encoded by $w_{tp}$, and constraints~\eqref{mrt:w} enforce the nonnegativity condition.
\begin{subequations}
    \begin{alignat}{2}
      &w_{1, p} = 1 - y_{1,p}, \quad &\forall p \in \mathcal{P}\label{mrt:w0}\\
      &w_{tp}\le T(1-y_{tp}), \quad & \forall t\in \mathcal{T}, \forall p \in \mathcal{P}\label{mrt:counter_1} \\
      &w_{tp} - w_{t-1,p} \le 1,\quad & \forall t\in \mathcal{T}\setminus\{1\}, \forall p \in \mathcal{P}\label{mrt:counter_2} \\
      &w_{tp} - w_{t-1,p} \ge 1- Ty_{tp}, \quad & \forall t\in \mathcal{T}\setminus\{1\}, \forall p \in \mathcal{P}\label{mrt:counter_3}\\
      &w_{tp}\le Z,\quad & \forall t\in \mathcal{T}, \forall p \in \mathcal{P}\label{mrt:z_w}\\
      &w_{tp}\ge0,\quad & \forall t\in \mathcal{T}, \forall p \in \mathcal{P}\label{mrt:w}
\end{alignat}
\end{subequations}

In summary, the MMRT formulation is given as:
\begin{alignat}{2}
    \textsc{(\hypertarget{mmrt}{MMRT})} \quad \min \quad &Z \label{mrt:obj}\\
    \text{s.t.}\quad &\sum_{j\in \mathcal{J}} x_{j} = N \notag \\
    & \sum_{j \in \mathcal{J}}V_{tjp}x_{j} - r_{tp}\ge -J(1-y_{tp}) - \beta,  &\quad \forall t \in \mathcal{T}, \forall p \in \mathcal{P}\notag\\
    & \sum_{j \in \mathcal{J}}V_{tjp}x_{j} - r_{tp}\le Jy_{tp} - \beta,  &\quad \forall t \in \mathcal{T}, \forall p \in \mathcal{P}\notag\\
    &w_{1, p} = 1 - y_{1,p}, &\quad\forall p \in \mathcal{P}\notag\\
    &w_{tp}\le T(1-y_{tp}),  &\quad \forall t\in \mathcal{T}, \forall p \in \mathcal{P}\notag \\
    &w_{tp} - w_{t-1,p} \le 1,  &\quad \forall t\in \mathcal{T}\setminus\{1\}, \forall p \in \mathcal{P}\notag \\
    &w_{tp} - w_{t-1,p} \ge 1- Ty_{tp},  &\quad \forall t\in \mathcal{T}\setminus\{1\}, \forall p \in \mathcal{P}\notag\\
    &w_{tp}\le Z, &\quad \forall t\in \mathcal{T}, \forall p \in \mathcal{P}\notag\\
    &w_{tp}\ge0, &\quad \forall t\in \mathcal{T}, \forall p \in \mathcal{P}\notag\\
   &x_{j} \in \{0, 1\}, &\quad \forall j\in \mathcal{J} \notag\\
    & y_{tp} \in \{0,1\}, &\quad  \forall t\in \mathcal{T}, \forall p \in\mathcal{P} \notag\\
    & Z \in\mathbb{Z}_{\ge0}
\end{alignat}

\begin{remark}[Cyclic Property]
    If the constellation coverage timeline $\bm{b}_{p}$ is cyclic, and its repetition period coincides with the mission horizon, as in the case of using common RGT orbital slots \cite{lee2020satellite}, then the first and last time steps of the optimization problem are coupled, yielding a cyclic visibility profile. To account for this property, constraints~\eqref{mrt:w0} are replaced by constraints~\eqref{mrt:circularity_1} and~\eqref{mrt:circularity_2} to link the first and last time steps.
    \begin{subequations}
        \begin{alignat}{2}
            &w_{1,p} - w_{T,p} \le 1,\quad & \forall p \in \mathcal{P}\label{mrt:circularity_1} \\
            &w_{1,p} - w_{T,p} \ge 1- Ty_{1,p}, \quad & \forall p \in \mathcal{P}\label{mrt:circularity_2}
        \end{alignat}
    \end{subequations}
\end{remark}

\begin{remark}[Minimizing the Sum of MRTs]
    If, instead of minimizing the MRT of the entire set of targets $\mathcal{P}$, the user seeks to minimize the sum of each target's MRT, the formulation has to be modified as follows. First, we drop decision variable $Z$ and propose decision variables $Z_{p} \in \mathbb{Z}_{\ge 0}$ to encode target $p$'s MRT. Then, objective function~\eqref{mrt:obj} is recast as objective function~\eqref{mrt:obj2}, and constraints~\eqref{mrt:z_w} are formulated as constraints~\eqref{mrt:z_w_2} to couple each target's $p$ decision variable $Z_{p}$ with the maximum coverage gap encoded in $w_{tp}$.
    \begin{alignat}{2}
    &\sum_{p \in \mathcal{P}}Z_{p}\label{mrt:obj2}\\
    &w_{tp}\le Z_{p},\quad & \forall t\in \mathcal{T}, \forall p \in \mathcal{P}\label{mrt:z_w_2}
\end{alignat}
\end{remark}

\begin{figure}[htbp]
	\centering
	\begin{subfigure}[c]{0.492\linewidth}
		\centering
		\includegraphics[width=\linewidth]{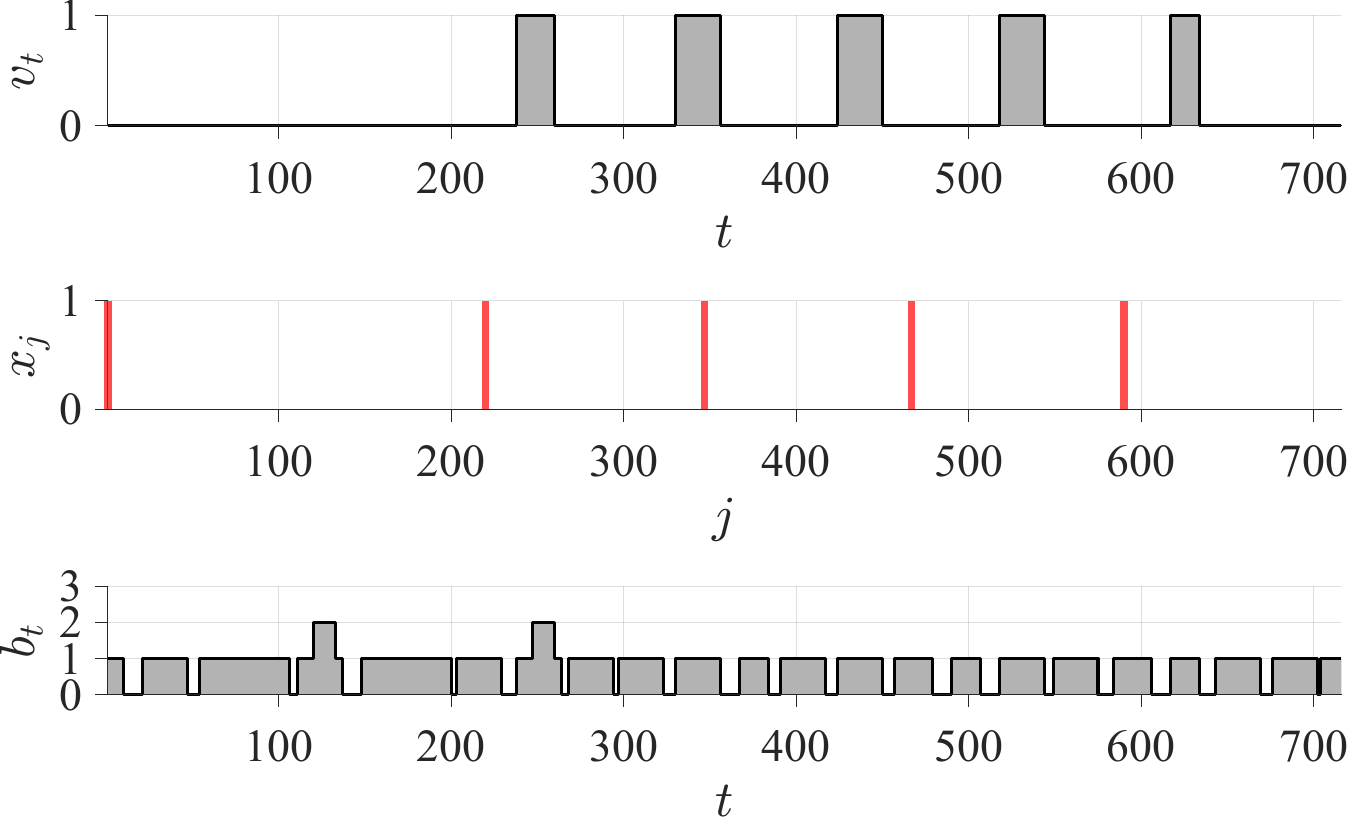}
		\caption{Illustration of the APC decomposition.}
            \label{fig:mrt_a}
	\end{subfigure}
	\begin{subfigure}[c]{0.492\linewidth}
		\centering
		\includegraphics[width=\linewidth]{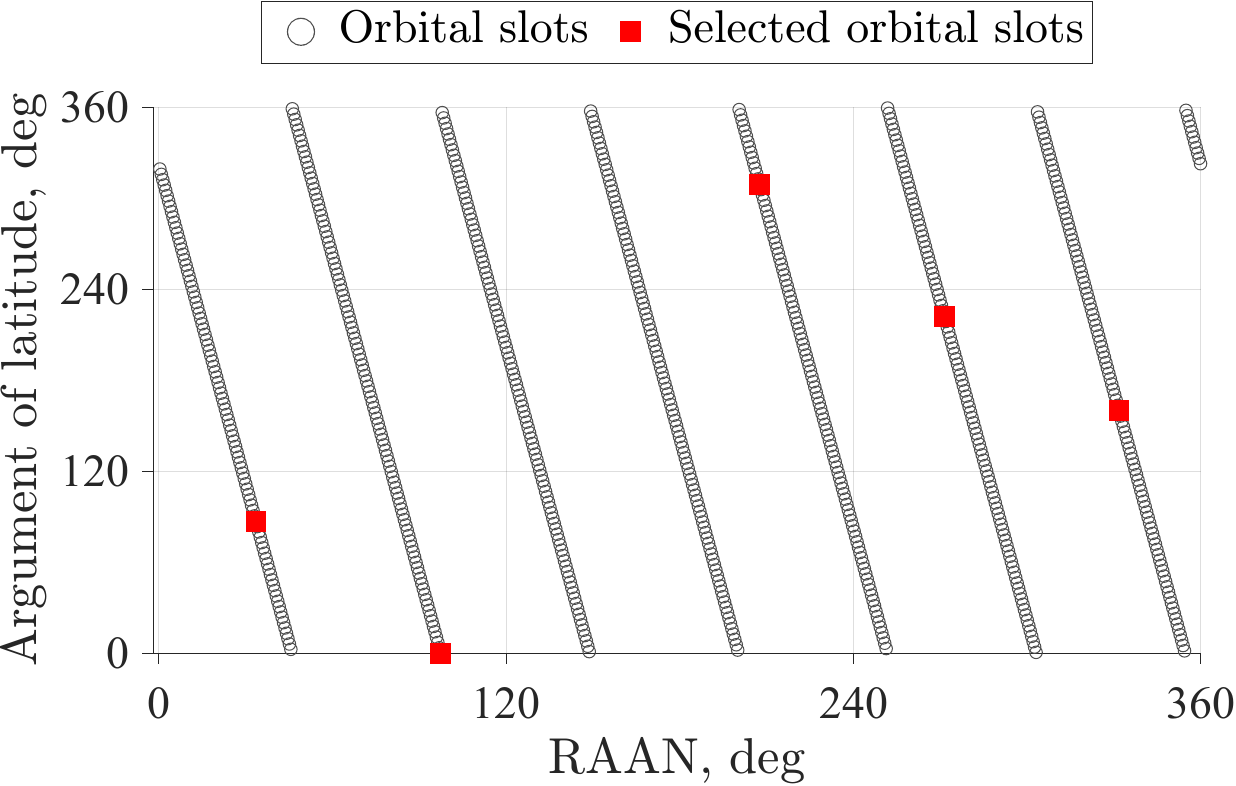}
		\caption{Optimal constellation configuration's orbital slot distribution.}
            \label{fig:mrt_b}
	\end{subfigure}
    \begin{subfigure}[c]{0.5853\linewidth}
            \centering
             \includegraphics[width=0.60\linewidth]{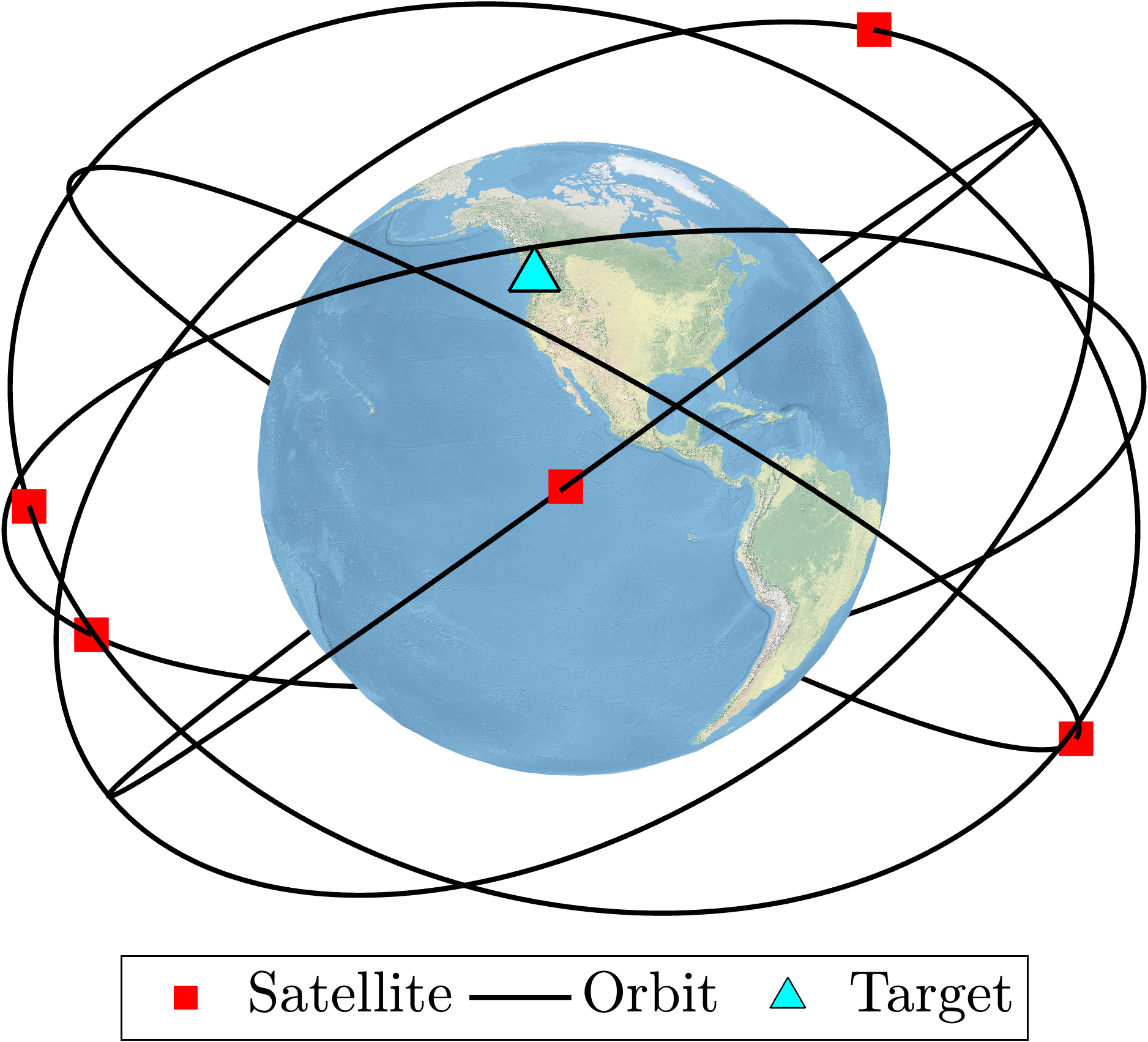}
            \caption{Five-satellite optimal constellation configuration with 22 min MRT.}
                \label{fig:mrt_c}
        \end{subfigure}
	\caption{MMRT formulation results.}
	\label{fig:mrt}
\end{figure}

\begin{example}[Five-Satellite MMRT]
Given a fixed number of satellites equal to five, the optimal objective value is 11, and the runtime is \SI{151.16}{s}, corresponding to an optimal constellation configuration with an MRT of \SI{22}{min} over Seattle. Figure~\ref{fig:mrt_a} presents the APC decomposition, and Fig.~\ref{fig:mrt_b} outlines the distribution of the orbital slots in the RAAN versus argument of latitude plane. Figure~\ref{fig:mrt_c} illustrates a 3D visualization of the constellation configuration and the target at the epoch.
\end{example}

\subsection{Minimal Average Revisit Time Problem}\label{sec:mart}
MART determines the location of $N$ satellites within a set $\mathcal{J}$ to minimize the ART. The value of MART lies in its ability to tackle the ART, a paramount metric in satellite constellation missions, where its success and performance depend on the average duration of the coverage gaps. The ART differs from MCLP and MMRT in that it focuses solely on the ART, rather than the collected observation rewards or the maximum duration of a coverage gap. Table~\ref{table:parameters_art} lists the decision variables specific to MART.

\begin{table}[htbp]
\caption{Specific MART decision variables.}
\centering
\begin{tabular}{lll}
\hline
\hline
Type & Symbol & Description \\
\hline
Decision variables  &  $\ell_{tp}$ & $\begin{cases}
                     1, &\text{if a coverage gap starts at time step $t$ for target $p$ } \\
                    0, &\text{otherwise}
                    \end{cases}$ \\
                    & $\alpha_{p}$ & Upper bound of target $p$'s ART \\
                    & $s_{tp}$ & Auxiliary variable for target $p$ at time step $t$ \\                    
\hline
\hline
\label{table:parameters_art}
\end{tabular}
\end{table}

To encode the start of a coverage gap during the mission, we propose gap indicator variables $\bm{\ell}:=\{\ell_{tp} \in \{0,1\}:t \in \mathcal{T},p \in \mathcal{P}\}$, where each element is defined as:
\begin{equation}
    \ell_{tp} = \begin{cases}
    1, & \text{if a coverage gap starts at time step $t$ for target $p$} \\
    0, & \text{otherwise}\\
\end{cases}
\end{equation}

The ART is expressed as the sum of all coverage gap durations divided by the number of gaps, defined as:
\begin{equation}
    \text{ART} :=\sum_{p\in\mathcal{P}}\frac{T-\sum_{t\in \mathcal{T}}y_{tp}}{\sum_{t\in \mathcal{T}} \ell_{tp}} \label{eq:art_definition}
\end{equation}
Note that this equation holds when the denominator is greater than zero. Given the nonlinear nature of the ART and considering that we aim to tackle the problem with an MILP formulation, we proceed to linearize it. First, we propose target $p$'s decision variable $\alpha_{p} \in \mathbb{R}_{\ge0}$, where its summation as an upper bound of the ART:
\begin{equation}
    \sum_{p \in \mathcal{P}}\alpha_{p} \ge\sum_{p\in\mathcal{P}}\frac{T-\sum_{t\in \mathcal{T}}y_{tp}}{\sum_{t\in \mathcal{T}} \ell_{tp}}\label{eq:art_ub}
\end{equation}
Then, we introduce auxiliary decision variables $s_{tp}\in\mathbb{R}_{\ge0}$ for target $p$ at time step $t$, defined as:
\begin{equation}
    s_{tp} = \alpha_{p}\ell_{tp}
\end{equation}
which couples the left-hand side and right-hand side's denominator of Eq.~\eqref{eq:art_ub}. To express this nonlinear equality as a set of linear inequality constraints, we leverage the Big-M method. Then, the new set of constraints that linearizes Eq.~\eqref{eq:art_ub} is:
\begin{subequations}
    \begin{alignat}{2}
        &s_{tp} \le T\ell_{tp}, \quad & \forall t \in \mathcal{T}, \forall p \in \mathcal{P}\label{eq:art_a_bm1}\\
        &s_{tp}\le \alpha_{p}, \quad & \forall t \in \mathcal{T}, \forall p \in \mathcal{P}\label{eq:art_a_bm2}\\
        &s_{tp}\ge \alpha_{p} - T(1-\ell_{tp}), \quad & \forall t \in \mathcal{T}, \forall p \in \mathcal{P}\label{eq:art_a_bm3}
    \end{alignat}
\end{subequations}
where the Big-M method's parameter takes the value of $T$. Subsequently, we couple for each target $p$ auxiliary variables $s_{tp}$ with the nominator of Eq.~\eqref{eq:art_definition} as:
\begin{equation}
    \sum_{t \in \mathcal{T}}s_{tp}\ge T-\sum_{t \in \mathcal{T}}y_{tp}, \quad\forall p \in \mathcal{P}\label{eq:art_a_y}
\end{equation}

Leveraging the linearization of Eq.~\eqref{eq:art_definition}, MART minimizes the ART with objective function:
\begin{equation}
    \sum_{p \in \mathcal{P}}\alpha_{p} \label{art:obj}
\end{equation}

To account for each target $p$'s coverage state, we leverage constraints~\eqref{mrt:bigM_1} and~\eqref{mrt:bigM_2}. In addition, we propose a set of constraints, based on the Big-M method, to couple gap indicator variables $\bm{\ell}$ with coverage state decision variables $\bm{y}$. If target $p$ is covered at time step $t$, constraints~\eqref{art:i_sum1} require gap indicator variables $\ell_{tp}=0$. Conversely, constraints~\eqref{art:i_sum2} enforce gap indicator variables $\ell_{tp}$ to be one if a coverage gap starts. Further, constraints~\eqref{art:i_sum3} enforce decision variables $\ell_{tp}$ to zero inside a coverage gap.
\begin{subequations}
    \begin{alignat}{2}
    & \ell_{tp} \le 1-y_{tp}, \quad & \forall t\in \mathcal{T}, \forall p \in\mathcal{P}\label{art:i_sum1}\\
    & \ell_{tp} \ge y_{t-1,p} - Ty_{tp}, \quad & \forall t\in \mathcal{T}\setminus\{1\}, \forall p \in\mathcal{P}\label{art:i_sum2}\\
    & \ell_{tp} \le y_{t-1,p} + y_{tp}, \quad & \forall t\in \mathcal{T}\setminus\{1\}, \forall p \in\mathcal{P}\label{art:i_sum3}
\end{alignat}
\end{subequations}

Lastly, piecing all constraints and the objective function together, the MART formulation is given as:
\begin{alignat}{2}
    \textsc{(\hypertarget{mart}{MART})} \quad \min \quad & \sum_{p\in\mathcal{P}}\alpha_{p} \notag\\
    \text{s.t.}\quad &\sum_{j\in \mathcal{J}} x_{j} = N \notag \\
    & \sum_{j \in \mathcal{J}}V_{tjp}x_{j} - r_{tp}\ge -J(1-y_{tp}) - \beta,  &\quad \forall t \in \mathcal{T}, \forall p \in \mathcal{P} \notag\\
    & \sum_{j \in \mathcal{J}}V_{tjp}x_{j} - r_{tp}\le Jy_{tp} - \beta,  &\quad \forall t \in \mathcal{T}, \forall p \in \mathcal{P}\notag\\
     &s_{tp} \le T\ell_{tp},  &\quad \forall t \in \mathcal{T}, \forall p \in \mathcal{P} \notag\\
    &s_{tp}\le \alpha_{p},  &\quad \forall t \in \mathcal{T}, \forall p \in \mathcal{P} \notag\\
    &s_{tp}\ge \alpha_{p} - T(1-\ell_{tp}),  &\quad \forall t \in \mathcal{T}, \forall p \in \mathcal{P} \notag\\
    &\sum_{t \in \mathcal{T}}s_{tp}\ge T-\sum_{t \in \mathcal{T}}y_{tp}, &\quad\forall p \in \mathcal{P} \notag\\
    & \ell_{tp} \le 1-y_{tp}, & \quad \forall t\in \mathcal{T}, \forall p \in\mathcal{P}\notag\\
    & \ell_{tp} \ge y_{t-1,p} - Ty_{tp}, &\quad \forall t\in \mathcal{T}\setminus\{1\}, \forall p \in\mathcal{P} \notag\\
    & \ell_{tp} \le y_{t-1,p} + y_{tp},  &\quad \forall t\in \mathcal{T}\setminus\{1\}, \forall p \in\mathcal{P} \notag\\
    & x_{j} \in \{0, 1\},  & \quad \forall j\in \mathcal{J} \notag\\
    & y_{tp} \in \{0,1\}, & \quad \forall t\in \mathcal{T}, \forall p \in\mathcal{P} \notag\\
    & \ell_{tp} \in \{0,1\}, &\quad\forall t\in \mathcal{T}, \forall p \in\mathcal{P}\label{art:i}\\
    & \alpha_{p} \ge 0, & \quad\forall p \in\mathcal{P} \label{art:alpha_p}\\
    & s_{tp} \ge 0,  &\quad\forall t\in \mathcal{T}, \forall p \in\mathcal{P} \label{art:s_tp} 
\end{alignat}

\begin{remark}[Cyclic Property]
To incorporate the periodicity of the constellation coverage timeline $\bm{b}_{p}$, constraints~\eqref{art:circularity_2} and~\eqref{art:circularity_3} are introduced to couple the first and last time steps.
\begin{subequations}
    \begin{alignat}{2}
        & \ell_{1,p} \ge y_{T,p}-Ty_{1,p}, \quad &\forall p \in\mathcal{P}\label{art:circularity_2}\\
        & \ell_{1,p} \le y_{T,p} + y_{1,p}, \quad &\forall p \in\mathcal{P}\label{art:circularity_3}
    \end{alignat}
\end{subequations}
\end{remark} 

\begin{figure}[htbp]
	\centering
	\begin{subfigure}[c]{0.492\linewidth}
		\centering
		\includegraphics[width=\linewidth]{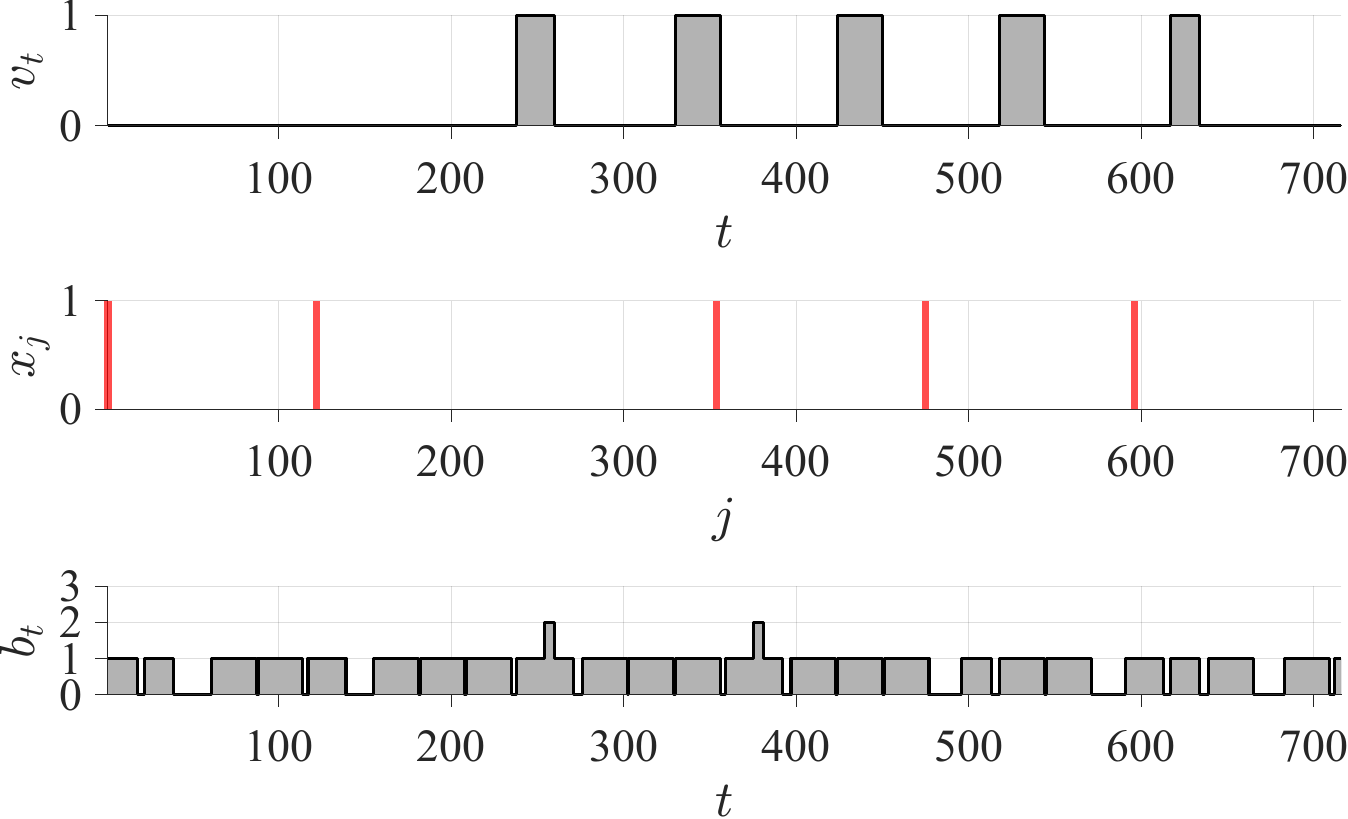}
		\caption{Illustration of the APC decomposition.}
            \label{fig:art_a}
	\end{subfigure}
	\begin{subfigure}[c]{0.492\linewidth}
		\centering
		\includegraphics[width=\linewidth]{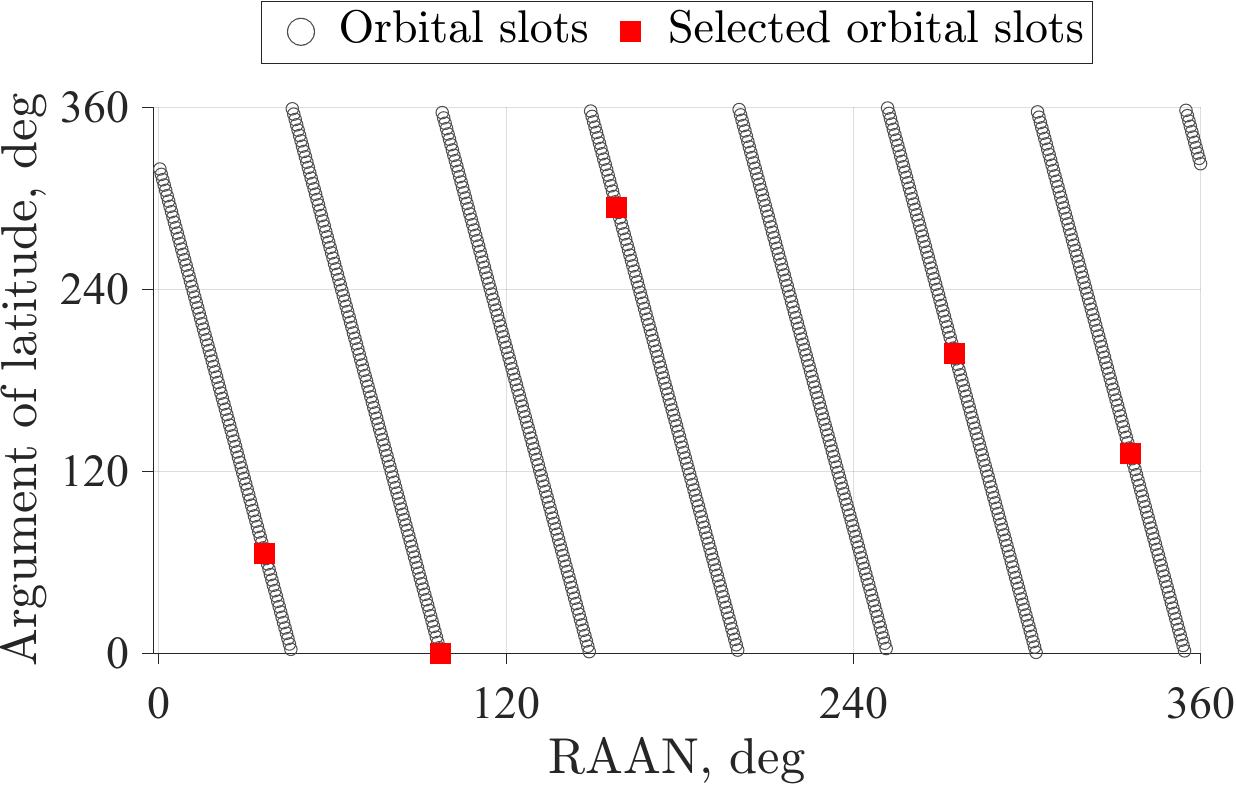}
		\caption{Optimal constellation configuration's orbital slot distribution.}
            \label{fig:art_b}
	\end{subfigure}
    \begin{subfigure}[c]{0.5853\linewidth}
            \centering
             \includegraphics[width=0.60\linewidth]{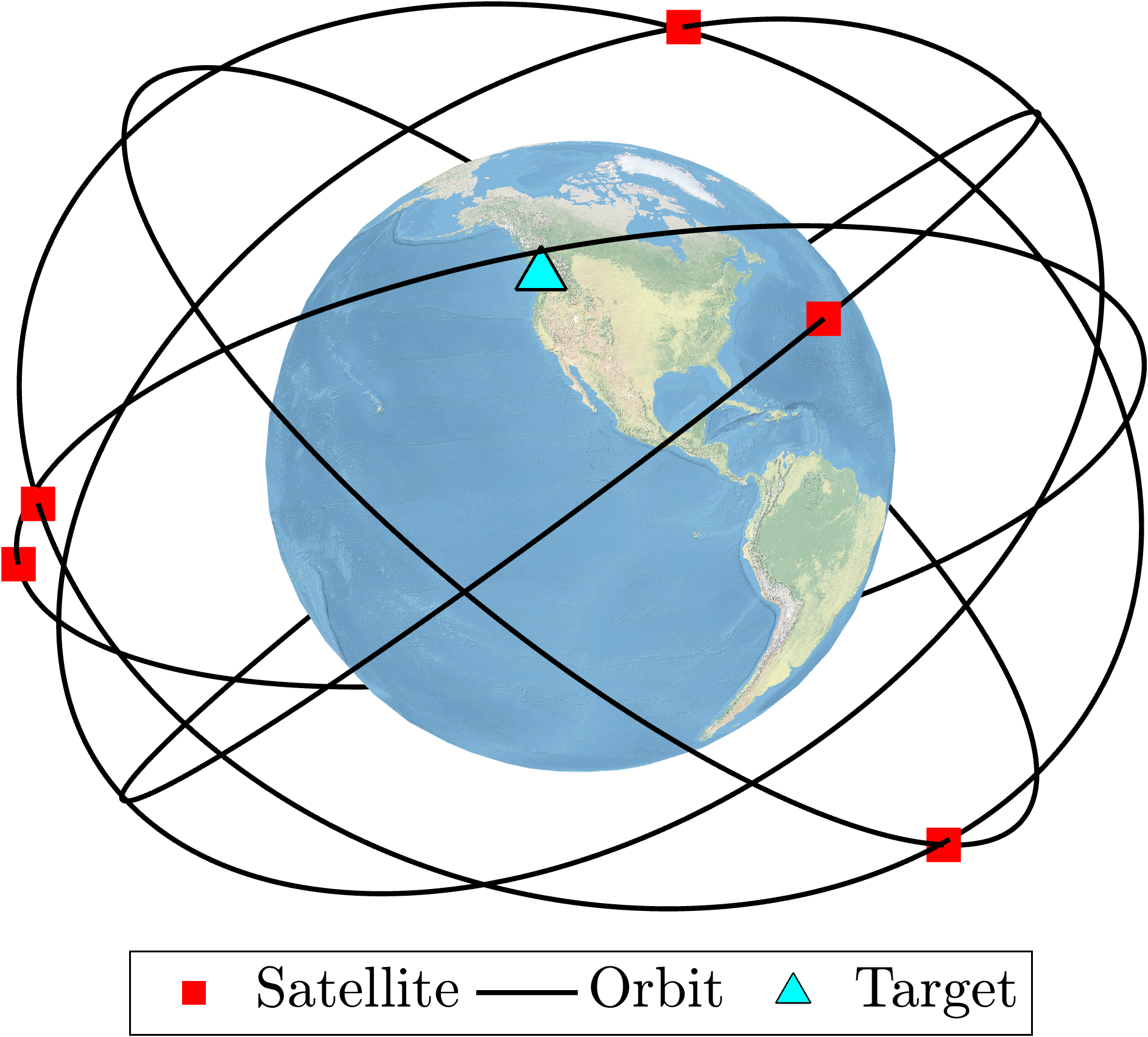}
            \caption{Five-satellite optimal constellation configuration with 12.43 min ART.}
                \label{fig:art_c}
        \end{subfigure}
	\caption{MART formulation results.}
	\label{fig:art}
\end{figure}

\begin{example}[Five-Satellite MART]
Given a fixed number of satellites equal to five, the optimal objective value is 6.21, and the runtime is \SI{819.14}{s}, corresponding to an optimal constellation configuration with an ART of \SI{12.43}{min} over Seattle. Figure~\ref{fig:art_a} presents the APC decomposition, and Fig.~\ref{fig:art_b} outlines the distribution of the orbital slots in the RAAN versus argument of latitude plane. Figure~\ref{fig:art_c} illustrates a 3D visualization of the constellation configuration and the target at the epoch.
\end{example}

\section{Comparative Analyses, Mix and Match, and Extensions}\label{sec:case_studies}
This section presents comparative analyses between the five formulations presented in Sec.~\ref{sec:mathematical_modeling} to demonstrate the differences between their optimization objectives and constraints, which lead to different constellation configurations and coverage metrics. Although the collection of MILP is useful in determining optimal constellation configurations tackling key coverage figures of merit, it may not be comprehensive for every mission scenario. Therefore, to account for specific mission settings, we tailor the MILP collection through a mix and match (\textit{i.e.,} interchanging the objective functions and constraints) and extensions (\textit{i.e., }incorporating new constraints or new parameter interpretations).

\subsection{Comparative Analyses of the Five Formulations}
This section performs two comparative analyses between the five formulations of Sec.~\ref{sec:mathematical_modeling}. First, we compare the number of satellites and percentage coverage for SCLP, PSCLP, and MCLP. The purpose of this comparative analysis is to highlight the impact of the constellation's number of satellites on the percentage coverage delivered. Second, we compare the discontinuous coverage figures of merit, that is, percent coverage, MRT, and ART, obtained by MCLP, MMRT, and MART. This comparison seeks to showcase that given a fixed number of satellites, their orbital distribution significantly impacts the discontinuous coverage metrics. Table~\ref{table:comp_review_v1} presents a review of the features considered in each formulation, discriminating them either as an optimization objective or constraints.

\begin{table}[htbp]
\caption{Optimization objectives and constraints summary for MILP formulations.}
\label{table:comp_review_v1}
\centering
\begin{threeparttable}
\begin{tabular}{llccccc}
\hline
\hline
Optimization & Figure of merit & \hyperlink{sclp}{SCLP} &\hyperlink{psclp}{PSCLP}& \hyperlink{mclp}{MCLP} & \hyperlink{mmrt}{MMRT} & \hyperlink{mart}{MART}\\
\hline
Objective & Minimize cost\tnote{*}  & \pmb{\checkmark}& \pmb{\checkmark}& & &\\
& Maximize observation reward\tnote{\textdagger} & & &\pmb{\checkmark} & &\\
& Minimize MRT & & & &\pmb{\checkmark} &\\
& Minimize ART & & & & &\pmb{\checkmark}\\
\hline
Constraints & Fixed number of satellites & & &\pmb{\checkmark} &\pmb{\checkmark} &\pmb{\checkmark}\\
& Spatiotemporal coverage requirements &\pmb{\checkmark} &\pmb{\checkmark} & & &\\
\hline
\hline
\end{tabular}
\begin{tablenotes}
\item[*] For example, the number of satellites, deployment cost (Sec.~\ref{sec:ext_deployment}).
\item[\textdagger] For example, percent coverage.
\end{tablenotes}
\end{threeparttable}
\end{table}

Table~\ref{table:comp_1} presents the number of satellites and the delivered coverage of the optimal constellations obtained by SCLP, PSCLP, and MCLP, where gray cells indicate the optimization's objective. The SCLP constellation requires eight satellites to deliver continuous coverage, while the PSCLP constellation requires six satellites to deliver at least \SI{85}{\%} coverage, indicating that the \SI{15}{\%} difference is achieved by adding two satellites to the constellation. The comparison between the MCLP and PSCLP constellations shows the trade-off between constellation configuration cost (i.e., number of satellites), and the delivered coverage.

\begin{table}[htbp]
\caption{Comparison between SCLP, PSCLP, and MCLP results.}
\label{table:comp_1}
\centering
\begin{tabular}{llll}
\hline
\hline
Metric &\hyperlink{sclp}{SCLP} &\hyperlink{psclp}{PSCLP}& \hyperlink{mclp}{MCLP}\\
\hline
Number of satellites & \cellcolor[HTML]{e0e0e0}8 & \cellcolor[HTML]{e0e0e0}6 & 5 (req. 5) \\
Percent coverage & \SI{100.0}{\%} (req. 100\%)  & \SI{85.75}{\%} (req. 85\%)  & \cellcolor[HTML]{e0e0e0}\SI{80.86}{\%}\\
\hline
\hline
\end{tabular}
\end{table}
 
Table~\ref{table:comp_2} presents the percentage coverage, MRT, and ART obtained by MCLP, MMRT, and MART. Analogously to the previous table, the gray cells indicate the objective value obtained by each formulation. It can be observed that each formulation outperforms the other two with respect to the metric considered in its corresponding objective function. MCLP provides the theoretical maximum percent coverage achievable by a constellation. Specifically, MCLP obtains \SI{80.86}{\%} percentage coverage while MMRT and MART have \SI{78.07}{\%} and \SI{80.02}{\%}, respectively. Similarly, MMRT and MART obtain the minimum MRT and ART, respectively. More precisely, MMRT has an MRT of \SI{22}{min}, while MCLP and MART have MRTs of \SI{52}{min} and \SI{44}{min}, respectively. Furthermore, MART has an ART of \SI{12.43}{min}, the lowest compared to the \SI{18.26}{min} and \SI{14.95}{min} obtained by the MCLP and MMRT formulations, respectively. The results presented in this table enable us to conclude that given a fixed number of satellites in the constellation, improving one of the discontinuous coverage's figures of merit does not imply an improvement in the other figures of merit. Thus indicating that there is necessarily no one configuration that is best in all metrics. For instance, it is clear how MMRT sacrifices percentage coverage by breaking the coverage timeline, aiming to reduce the duration of the coverage gaps.

\begin{table}[htbp]
\caption{Comparison between MCLP, MMRT, and MART results.}
\label{table:comp_2}
\centering
\begin{tabular}{llll}
\hline
\hline
Metric & \hyperlink{mclp}{MCLP} & \hyperlink{mmrt}{MMRT} & \hyperlink{mart}{MART}\\
\hline
Percent coverage & \cellcolor[HTML]{e0e0e0}\SI{80.86}{\%} & \SI{78.07}{\%}             & \SI{80.02}{\%}\\
MRT              & \SI{52}{min}                           & \cellcolor[HTML]{e0e0e0}\SI{22}{min}  & \SI{44}{min} \\
ART              & \SI{18.26}{min}                         & \SI{14.95}{min}        &\cellcolor[HTML]{e0e0e0}\SI{12.43}{min}  \\
\hline
\hline
\end{tabular}
\end{table}

\subsection{Mix and Match}\label{sec:add_ons}
This section intends to provide additional tools to the user by extending the applicability of the formulations presented in Sec.~\ref{sec:mathematical_modeling}. To this end, we conduct a mix and match between the constraints and the objective functions of SCLP, MMRT, and MART. Each formulation is accompanied by an illustrative example adopting the parameters defined in Sec.~\ref{sec:parameters}. 

\subsubsection{Maximum Revisit Time as a Mission Requirement}
We propose a mathematical formulation to design a minimum-cost constellation configuration where the MRT is a mission constraint. The new formulation, denoted as Constrained MRT (CMRT), adopts SCLP's objective function~\eqref{sclp:obj}, and leverages MMRT's constraints except for constraints~\eqref{mclp:cardinality}, which are dropped. To enforce the MRT as a mission requirement, decision variable $Z$ is recast as parameter $z_{p} \in \mathbb{Z}_{\ge0}$, which represents the upper bound on target $p$'s MRT. Consequently, constraints~\eqref{mrt:z_w} are replaced by constraints~\eqref{eq:sclp_mrt_zw}. The full formulation is given as:
\begin{alignat}{2}
    \textsc{(CMRT)} \quad \min \quad & \text{Objective function \eqref{sclp:obj}} \notag\\
    \text{s.t.} \quad & \text{Constraints~\eqref{sclp:xj}, \eqref{psclp:ytp}, \eqref{mrt:bigM_1}, \eqref{mrt:bigM_2},}\notag\\
    &\text{\eqref{mrt:w0}, \eqref{mrt:counter_1}, \eqref{mrt:counter_2}, \eqref{mrt:counter_3}, \eqref{mrt:w}}\notag\\
    &w_{tp}\le z_{p},\quad  \forall t\in \mathcal{T}, \forall p \in \mathcal{P}\label{eq:sclp_mrt_zw}
\end{alignat}
It is noteworthy to mention that this formulation minimizes the constellation configuration cost given a user-defined MRT upper bound, which is distinct from directly minimizing the MRT. Therefore, the obtained solution is not guaranteed to have a minimum MRT.

\begin{remark}(Cyclic Property)
    If the constellation coverage timeline $\bm{b}_{p}$ is cyclic, constraints~\eqref{mrt:w0} are dropped, and constraints~\eqref{mrt:circularity_1} and~\eqref{mrt:circularity_2} are used.
\end{remark}

\begin{example}[10-min CMRT]
    Given a \SI{10}{min} MRT upper bound, the optimal objective value is six, and the runtime is \SI{288.26}{s}, corresponding to an optimal constellation configuration of six satellites. To demonstrate the differences between CMRT and MMRT, we summarize in Table~\ref{table:mam_mrt} the obtained results corresponding to the illustrative examples, where gray cells indicate objective values. Conversely to CMRT, MMRT determines an optimal constellation configuration, assuming the number of satellites fixed to five, that obtains a minimum MRT of \SI{22}{min}. Further, the results show that adding one more satellite into the constellation reduces the MRT by approximately \SI{14}{min}. In essence, both formulations enable the user to design optimal constellation configurations adopting MRT as the discontinuous coverage figure of metric. However, they differ in their objective; while the MMRT minimizes the MRT, the CMRT minimizes the constellation configuration cost.
    \begin{table}[htbp]
        \caption{Comparison between CMRT and MMRT results.}
        \label{table:mam_mrt}
        \centering
        \begin{tabular}{lll}
        \hline
        \hline
        Metric & CMRT &\hyperlink{mmrt}{MMRT}\\
        \hline
        Number of satellites & \cellcolor[HTML]{e0e0e0}6 & 5 (req. 5) \\
        MRT & \SI{8}{min} (req. $\le$ \SI{10}{min})  & \cellcolor[HTML]{e0e0e0}\SI{22}{min}\\
        \hline
        \hline
        \end{tabular}
        \end{table}
\end{example}

\subsubsection{Average Revisit Time as a Mission Requirement}
We propose the Constrained ART (CART) to design a minimum-cost constellation configuration such that the ART is lower than a user-defined requirement. Objective function~\eqref{sclp:obj} is adopted to minimize the cost of the constellation. From the set of constraints used in Sec.~\ref{sec:mart} to minimize the ART, we drop cardinality constraints~\eqref{mclp:cardinality} and recast objective function~\eqref{art:obj} as constraints~\eqref{eq:sclp_art} where parameter $\gamma_{p} \in \mathbb{R}_{\ge0}$ indicates the upper bound of the ART for target $p$. The full formulation is given as:
\begin{alignat}{2}
    \textsc{(CART)} \quad \min \quad & \text{Objective function \eqref{sclp:obj}}  \notag\\
    \text{s.t.} \quad & \text{Constraints~\eqref{sclp:xj}, \eqref{psclp:ytp}, \eqref{mrt:bigM_1}, \eqref{mrt:bigM_2}, \eqref{eq:art_a_bm1}, \eqref{eq:art_a_bm2}, \eqref{eq:art_a_bm3}}, \notag \\
    &\text{\eqref{eq:art_a_y}, \eqref{art:i_sum1}, \eqref{art:i_sum2}, \eqref{art:i_sum3}, \eqref{art:i}, \eqref{art:alpha_p}, \eqref{art:s_tp}} \notag \\
    & \alpha_{p} \le \gamma_{p}, \quad \forall p \in \mathcal{P} \label{eq:sclp_art}
\end{alignat}
It should be noted that this formulation minimizes the cost of the constellation configuration given a set of constraints, and does not minimize the ART. Therefore, the delivered ART is not guaranteed to be minimum.

\begin{remark}(Cyclic Property)
    If the constellation coverage timeline $\bm{b}_{p}$ is cyclic, constraints~\eqref{art:circularity_2} and~\eqref{art:circularity_3} are incorporated.
\end{remark}

\begin{example}[8-min CART]
    Given an \SI{8}{min} ART upper bound, the optimal objective value is six, and the runtime is \SI{462.43}{s}, corresponding to an optimal constellation configuration with six satellites. We present in Table~\ref{table:mam_art} the comparison between the results obtained for this illustrative example and for the corresponding MART one. Conversely to CART, given a requirement of fixing the number of satellites to five, MART determines a constellation configuration that has an ART of \SI{12.43}{min}, showcasing that removing one satellite leads to an increase of nearly \SI{4}{min}. In essence, both formulations are useful tools for mission designers tackling the problem of constellation configuration design, where the ART is the discontinuous coverage figure of merit.
    \begin{table}[htbp]
        \caption{Comparison between CART and MART results.}
        \label{table:mam_art}
        \centering
        \begin{tabular}{lll}
        \hline
        \hline
        Metric & CART &\hyperlink{mart}{MART}\\
        \hline
        Number of satellites & \cellcolor[HTML]{e0e0e0}6 & 5 (req. 5) \\
        ART & \SI{7.88}{min} (req. $\le$ \SI{8}{min})  & \cellcolor[HTML]{e0e0e0}\SI{12.43}{min}\\
        \hline
        \hline
        \end{tabular}
        \end{table}
\end{example}

\subsection{Extensions}
We present a set of extensions to showcase additional applications of the five MILP formulations. Although not comprehensive enough to span all possible applications, we believe that they provide the reader with good insight into common constellation configuration design problems. First, through a new parameter interpretation, we demonstrate the ability of the formulations to consider two families of orbital slots and static and dynamic targets with time-varying coverage requirements. Second, we illustrate the sensitivity of the constellation configuration design problem to distinct orbital slots' cost interpretations. Lastly, we propose a new set of constraints that enable the mission designer to incorporate ISL into the constellation configuration design optimization.

\subsubsection{Constellation Configuration Design for Space-Based Space Situational Awareness}
This case study tackles the problem of designing an optimal space-based SSA constellation configuration with multi-layer, heterogeneous orbital slot families, considering static and dynamic targets with time-varying coverage requirements. We adopt PSCLP to design a constellation to cover the International Space Station (ISS), and the three Deep Space Network (DSN) stations located in California, USA (geodetic coordinates \{$35.42\degree$N, $116.89\degree$W, \SI{0}{m}\}), Madrid, Spain (geodetic coordinates \{$40.43\degree$N, $4.24\degree$W, \SI{0}{m}\}), and Canberra, Australia (geodetic coordinates \{$35.40\degree$S, $148.98\degree$E, \SI{0}{m}\}). Given the mission epoch defined as January 1, 2025, 12:00:00.000 UTC, the ISS's state vector is determined using MATLAB's built-in function \texttt{propagateOrbit} \cite{MATLAB}, and the TLE information required as input for the function is obtained from CelesTrak \cite{celestrak}. We require continuous coverage for the ISS and the DSN station located in Canberra, and \SI{60}{\%} coverage for the remaining two DSN stations.

The first family of orbital slots corresponds to the lower layer, circular RGT orbits with an inclination of \SI{110}{deg} and 11:1 resonance ratio, corresponding to a semi-major axis of \SI{8541.36}{km} and a repetition period of \SI{86457.78}{s}. The second family of orbital slots corresponds to the upper layer; circular non-RGT orbits have a semi-major axis of \SI{8741.36}{km}, an inclination of \SI{40}{deg}, and their RAAN and argument of latitude uniformly discretized over the entire plane with 20 and 15 steps, respectively. The time horizon is set coincident with the repetition period, and the time step size is defined as \SI{120}{s}, leading to 720 time steps. For the ISS, and DSN stations located in California and Madrid, we define a coverage threshold parameter $r_{tp}=1$ for all $t \in \mathcal{T}$, that is, single-fold coverage. Conversely, for the DSN station located in Canberra, we define $r_{tp}=3$ for all $t\in [280,400]$, $r_{tp}=2$ for all $t\in [500,600]$, and $r_{tp}=1$ for all $t\in [1,280)\bigcup (600,720]$, that is, time-varying coverage requirements.

\begin{figure}[htbp]
	\centering
	\begin{subfigure}[H]{0.492\linewidth}
    		\centering
    		\includegraphics[width=\linewidth]{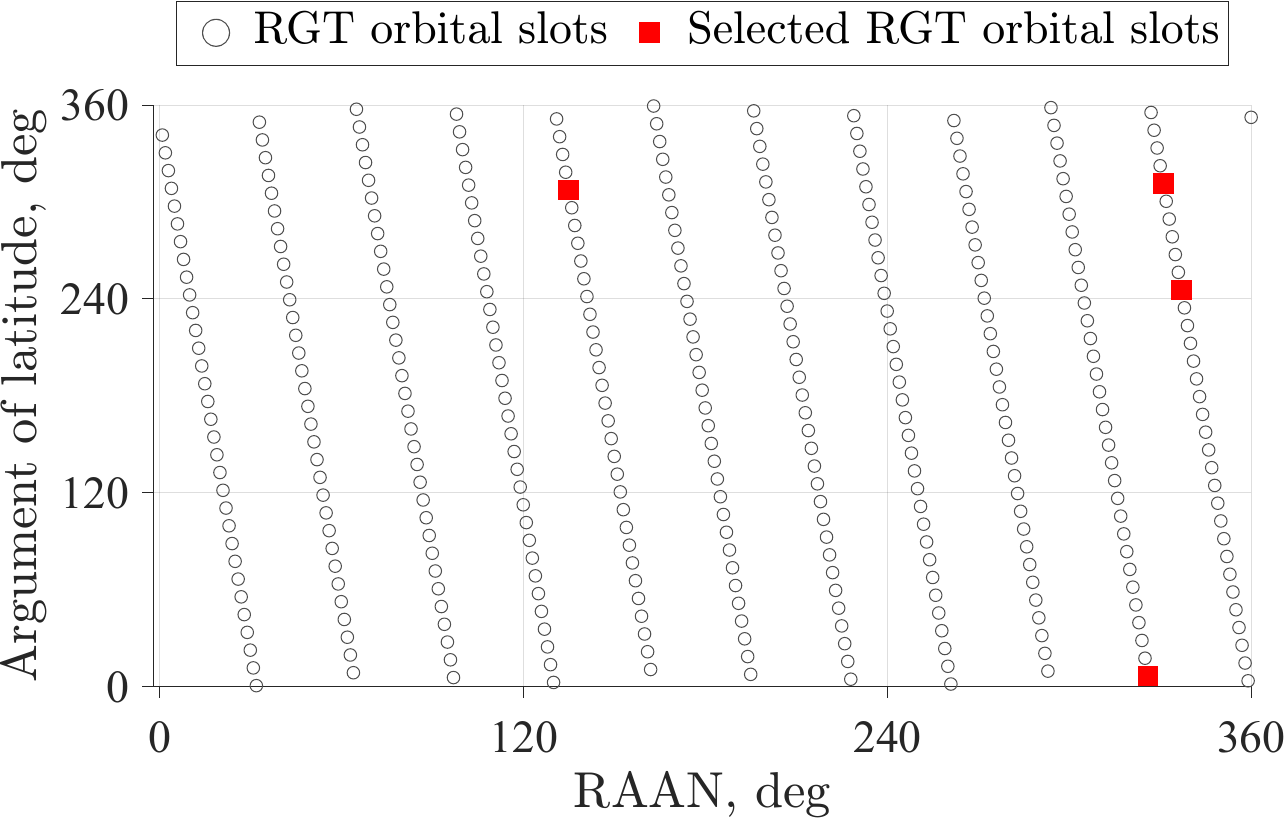}
    		\caption{Lower layer, RGT orbital slots distribution.}
                \label{fig:tracking_rgt}
	\end{subfigure}
    \begin{subfigure}[H]{0.492\linewidth}
    		\centering
    		\includegraphics[width=\linewidth]{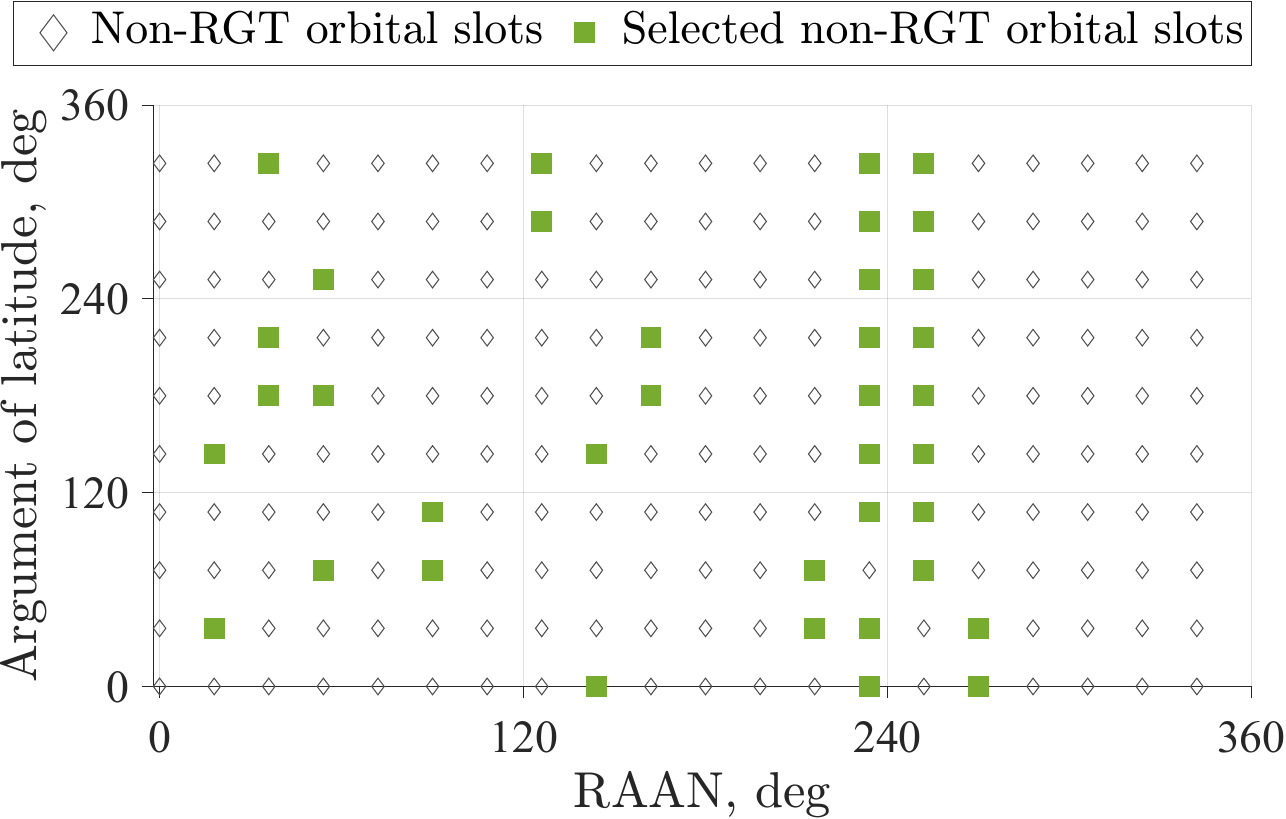}
    		\caption{Upper layer, non-RGT orbital slots distribution.}
                \label{fig:tracking_nonrgt}
	\end{subfigure}
	\caption{Optimal constellation configuration's orbital slot distribution.}
    \label{fig:tracking_pat_orb}
\end{figure}

\begin{figure}[htbp]
    \centering
    \includegraphics[width=0.9\linewidth]{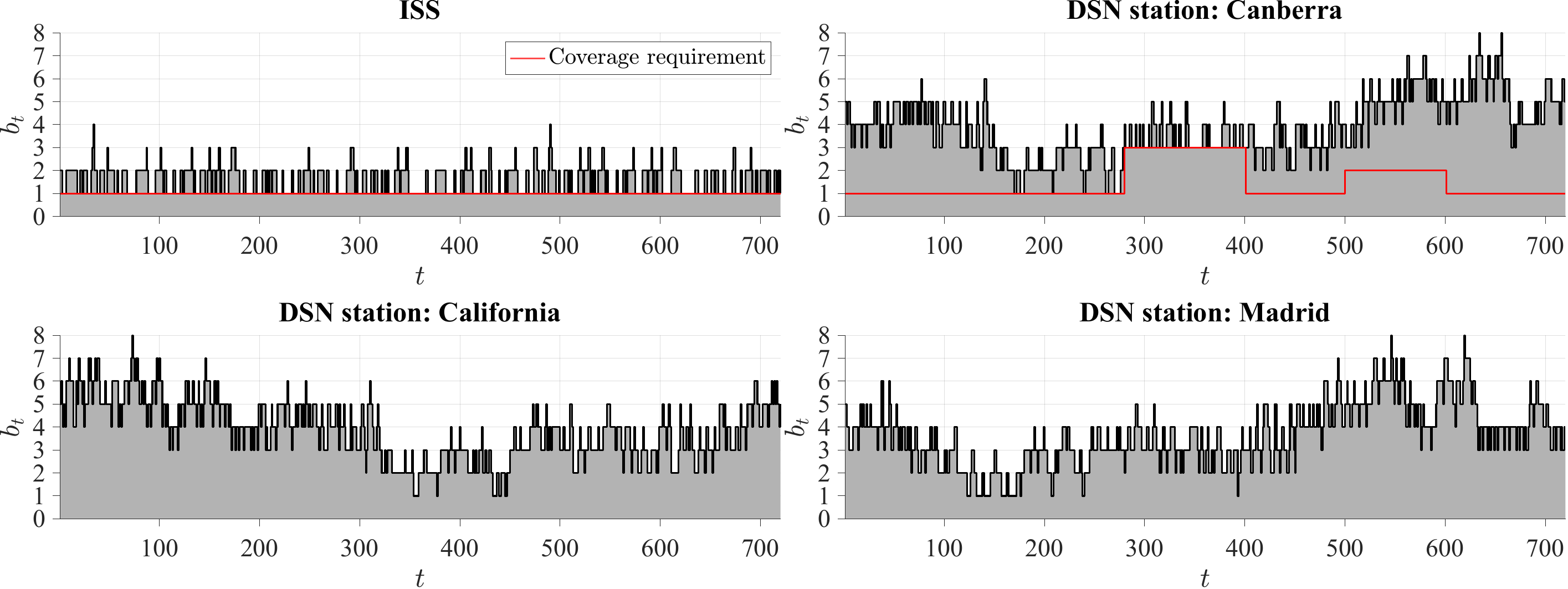}
    \caption{Coverage timeline over the targets.}
    \label{fig:tracking_cov}
\end{figure}

\begin{figure}[H]
    \centering
    \includegraphics[width=0.3\linewidth]{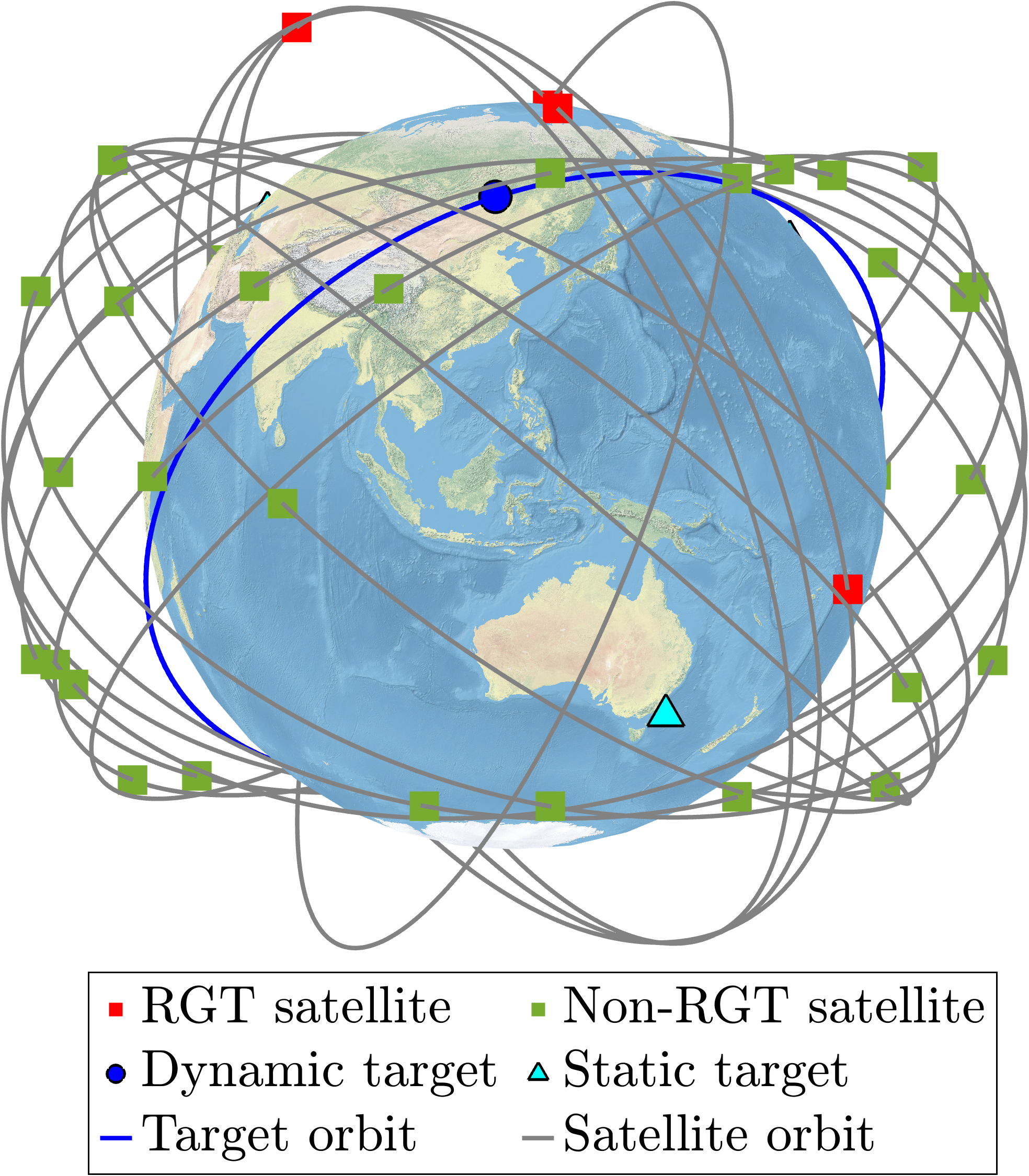}
    \caption{41-satellite optimal SSA constellation configuration.}
    \label{fig:tracking_3d}
\end{figure}

The optimal objective value is 41, and the runtime is \SI{3.92}{s}, corresponding to an optimal constellation of 41 satellites. From the total, four are RGT and 37 are non-RGT, where Fig.~\ref{fig:tracking_pat_orb} presents their distribution in the RAAN versus argument of latitude plane. From the figure, it is clear how a set of occupied non-RGT orbital slots is concentrated over the RAAN plane of \SI{240}{deg}, which is close to the ISS's RAAN of \SI{237.33}{deg}. Then, the remainder of the occupied orbital slots spans most of the RAAN versus argument of latitude plane. Figure~\ref{fig:tracking_cov} outlines the constellation coverage timeline $\bm{b}_{p}$ over the four targets, and the coverage threshold parameter $r_{tp}$, as a red line, for the targets with a continuous coverage requirement. It should be noted that although the DSN stations in California and Madrid did not impose a continuous coverage requirement, the optimal constellation delivers at least single-fold continuous coverage over them. Lastly, Fig.~\ref{fig:tracking_3d} showcases the optimal constellation configuration, and the targets with respect to the epoch.

We proceed to solve the problem by considering a single family of orbital slots and assuming the same set of parameters as described above. For the RGT orbital slots only, the problem is infeasible, which indicates that from this set is not possible to obtain a constellation configuration capable of satisfying the coverage requirements. For the non-RGT orbital slots only, the optimal constellation configuration uses 43 satellites, two more than for the case where the two families are considered. Therefore, this comparison provides information on the benefits of using heterogeneous families of orbital slots for the constellation configuration design optimization, as they complement each other.

\subsubsection{SCLP with Deployment Cost} \label{sec:ext_deployment}
This case study provides insight into the impact that different interpretations of each orbital slot's cost $c_{j}$ have on the optimal constellation configuration. First, the cost is interpreted as the $\Delta V$ required to transfer from the parking orbit to orbital slot $j$. Then, we define $c_{j}=1$ for all orbital slots $j \in \mathcal{J}$. We adopt the SCLP to design a minimum-cost constellation configuration to achieve single-fold continuous coverage over Morgantown, West Virginia, USA (geodetic coordinates \{$39.38\degree$N, $79.57\degree$W, \SI{0}{m}\}), and Seoul, Republic of Korea (geodetic coordinates \{$37.55\degree$N, $126.99\degree$E, \SI{0}{m}\}). Accordingly, the coverage threshold parameter is set as $r_{tp}=1$ for all $t \in \mathcal{T}$ and $p \in \mathcal{P}$, and a minimum elevation angle of \SI{10}{deg} is required for both targets. The epoch is designated as January 1, 2025, 12:00:00.000 UTC, the time step size is \SI{120}{s}, and the number of time steps is 1600.
 
The circular parking orbit is defined with a semi-major axis of \SI{6878.14}{km}, an inclination of \SI{28.5}{deg}, a RAAN of \SI{0}{deg}, and an argument of latitude of \SI{210}{deg}. For the orbital slots, 154 orbital planes with a semi-major axis of \SI{8378.14}{km} are defined from all possible combinations between the candidate RAAN and inclinations. The candidate RAANs are obtained by uniformly discretizing the range between 0 and \SI{360}{deg} in 14 steps. Similarly, the candidate inclinations encompass the range between 28.5 and \SI{37.55}{deg}, discretized in 11 intervals. In addition, each orbital plane contains 10 uniformly phased orbital slots.

To determine the cost of each orbital slot $j$, it is assumed that the orbital transfer vehicles (OTV) perform a Hohmman transfer to achieve a semi-major axis of \SI{8378.14}{km}. Subsequently, a change of plane maneuver is performed to arrive in the orbital plane with the corresponding RAAN and inclination. Lastly, phasing maneuvers are performed to cover the orbital plane. It is important to highlight that this case study aims to demonstrate the applicability of SCLP to obtain an optimal constellation configuration given specific transfer costs; consequently, a more efficient concept of operations for the deployment sequence can be applied.

\begin{figure}[htbp]
	\centering
	\begin{subfigure}[H]{0.492\linewidth}
    		\centering
    		\includegraphics[width=\linewidth]{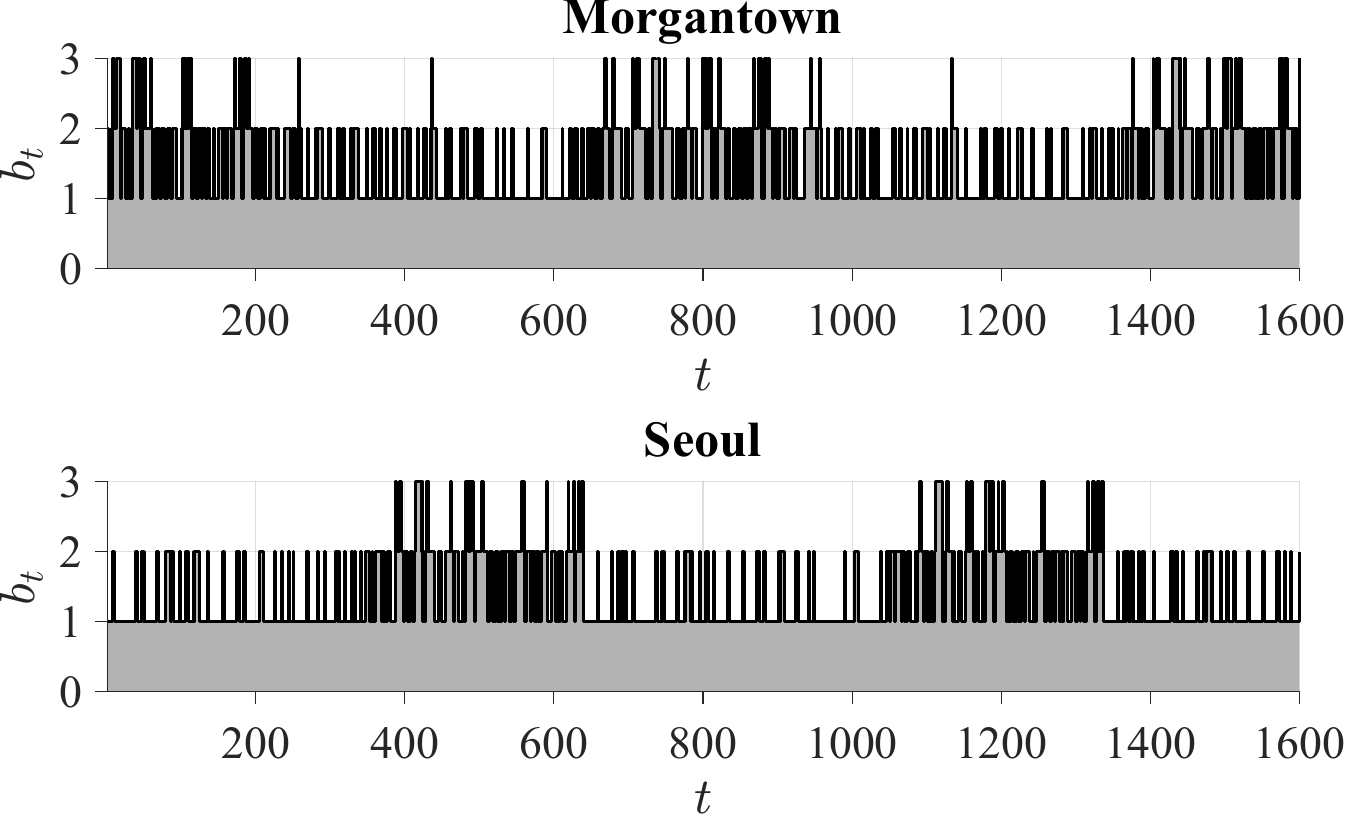}
    		\caption{Constellation pattern vector and coverage timeline.}
                \label{fig:deployment_a}
	\end{subfigure}
    \begin{subfigure}[H]{0.492\linewidth}
        \centering
        \includegraphics[width=0.56\linewidth]{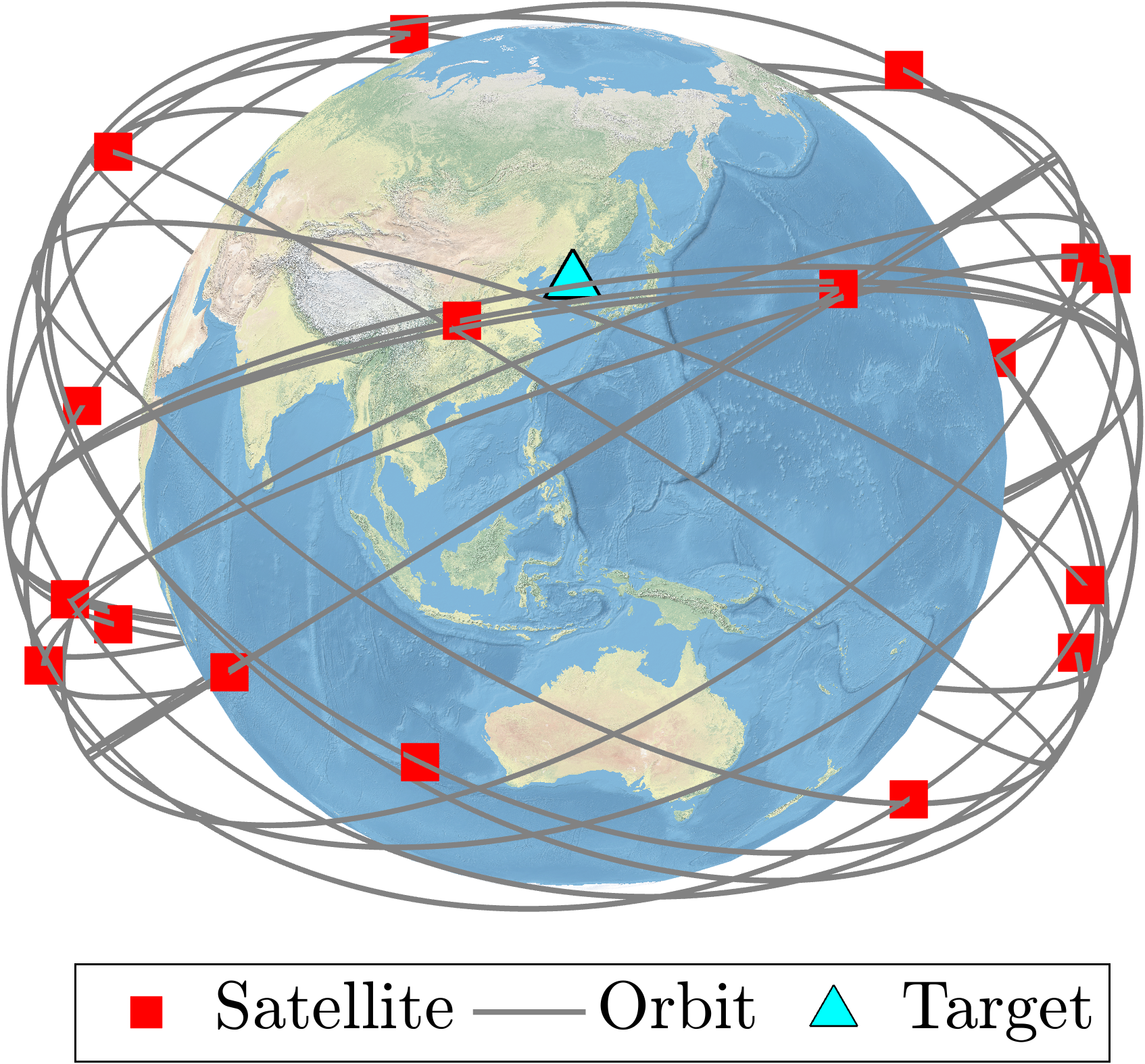}
        \caption{Minimum $\Delta V$ 20-satellite optimal constellation configuration.}
            \label{fig:deployment_e}
    \end{subfigure}
	\caption{SCLP with deployment cost results.}
	\label{fig:deployment_1}
\end{figure}

The optimal constellation configuration is obtained after \SI{112.27}{s}, has a cost of \SI{89.75}{km/s}, and uses 20 satellites. Figure~\ref{fig:deployment_a} outlines the constellation pattern vector $\bm{x}$ and the coverage timeline $\bm{b}_{p}$, delivering up to triple-fold coverage over both targets. Lastly, Fig.~\ref{fig:deployment_e} presents the 3D illustration of the minimum-cost constellation configuration.

\begin{figure}[htbp]
	\centering
	\begin{subfigure}[H]{0.46\linewidth}
    		\centering
    		\includegraphics[width=\linewidth]{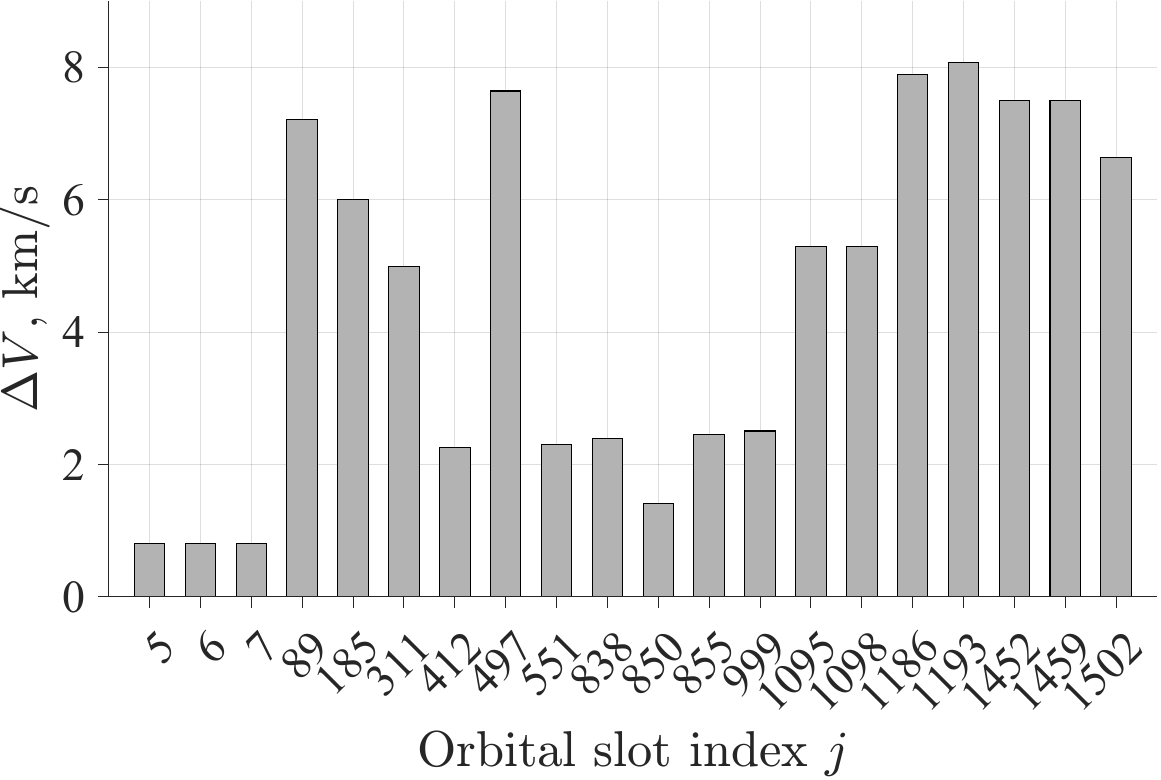}
    		\caption{20-satellite $\Delta V$ aware constellation with a cost of 89.75 km/s.}
                \label{fig:deployment_cost_aware}
	\end{subfigure}
    \begin{subfigure}[H]{0.46\linewidth}
        \centering
        \includegraphics[width=\linewidth]{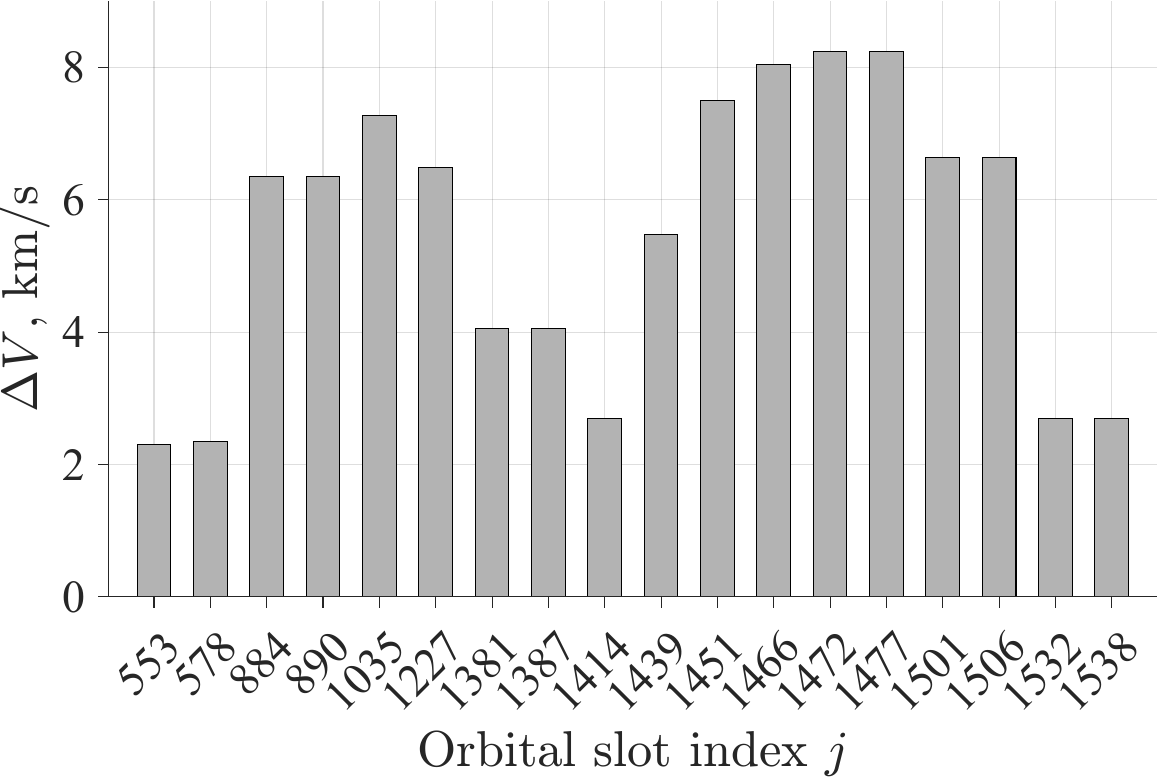}
        \caption{18-satellite $\Delta V$ non-aware constellation with a cost of 98.05 km/s.}
            \label{fig:deployment_cost_nonaware}
    \end{subfigure}
    \caption{Cost of occupied orbital slots for optimal configurations.}
        \label{fig:deployment_cost}
    \centering
\end{figure}

Defining a uniform cost for all orbital slots yields the constellation configuration that requires the minimum number of satellites to deliver continuous coverage over the targets. SCLP retrieves, after \SI{22.39}{s}, an optimal configuration that uses 18 satellites, two less than the optimal configuration that minimizes the deployment cost. The total $\Delta V$ required for the OTV to deploy the new constellation is \SI{98.05}{km/s}, representing an increase of \SI{9.24}{\%} with respect to the minimum $\Delta V$ required of \SI{89.75}{km/s}. Figure~\ref{fig:deployment_cost} presents the side-by-side comparison between the solution with minimum $\Delta V$, denoted as $\Delta V$ aware, and the solution with uniform cost, denoted as $\Delta V$ non-aware. The figure highlights how the non-aware constellation selects fewer orbital slots, but with larger $\Delta V$ values, leading to a higher total $\Delta V$ deployment cost.

\subsubsection{Inter-Satellite Links}
This section introduces a new set of constraints that enable the user to design an optimal constellation configuration with robust ISL. To model the link availability between two satellites, we propose the Boolean ISL parameter $W_{tjk}$, where each element is defined as:
\begin{equation}
W_{tjk}=\begin{cases}
1, &\text{if satellite $j$ can establish a link with satellite $k$ at time step $t$}\\
0, &\text{otherwise}
\end{cases} \label{isl:w}
\end{equation}
where in cases of mutual visibility, such as with bidirectional ISLs, $W_{tjk}=W_{tkj}$. 

\begin{figure}[htbp]
	\centering
		\includegraphics[width=0.3\linewidth]{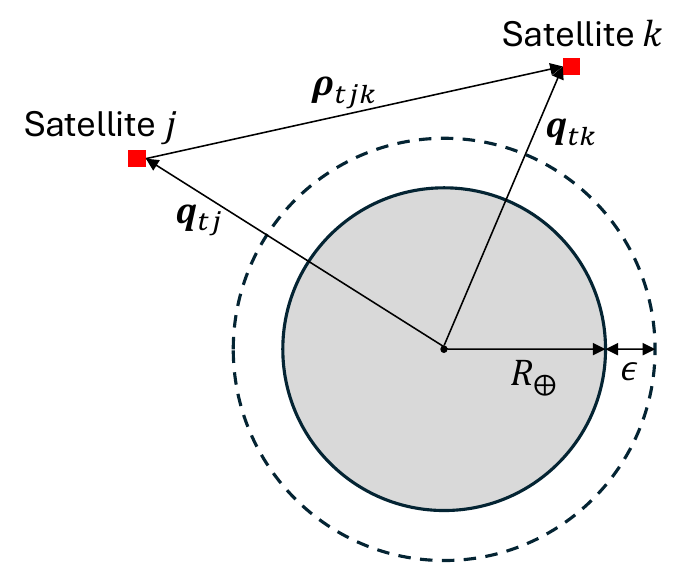}
		\caption{Geometrical relationship for inter-satellite visibility.}
		\label{fig:intersat}
\end{figure}

In this paper, to obtain $W_{tjk}$, we denote $\bm{q}_{tj}$ and $\bm{q}_{tk}$ as the position vectors of satellites $j$ and $k$ at time step $t$, respectively, and $\bm{\rho}_{tjk}$ denotes the relative position vector pointing from satellite $j$ to satellite $k$ at time step $t$. From the argument of geometry as depicted in Fig.~\ref{fig:intersat}, the inter-satellite visibility function $g_{tjk}$ is derived as follows:
\begin{equation}
    g_{tjk} = \big(q_{tj}^{2}-(R_{\oplus}+\epsilon)^{2}\big)^{1/2}+\big(q_{tk}^{2}-(R_{\oplus}+\epsilon)^{2}\big)^{1/2}-\rho_{tjk}
\end{equation}
where $R_{\oplus}=\SI{6378.14}{km}$ is the radius of the Earth, and $\epsilon$ is the bias factor to accommodate the additional margin of altitude (\textit{e.g.}, signal attenuation due to atmosphere). It is intuitive that when $g_{tjk}>0$, satellite $j$ has visibility over satellite $k$ at time step $t$. When $g_{jk}=0$, $\bm{\rho}_{tjk}$ is tangential to the sphere defined by the radius of $R_{\oplus}+\epsilon$. 

To generate an ISL network topology in which all satellites are connected, we propose leveraging Dirac's theorem \cite{Dirac} to guarantee the existence of a Hamiltonian cycle in the ISL network. Note that a cycle is Hamiltonian if it visits all vertices in the graph. To describe Dirac's theorem, we first define, for this section only, the following terms: $\mathcal{V}(G)$ represents the set of vertices on a graph $G$, $n$ is the number of vertices, and $d(\nu)$ is the degree of vertex $\nu \in \mathcal{V}(G)$.

\begin{theorem}[Dirac's Theorem \cite{Dirac}]
    Let $G$ be a simple graph with $n \ge 3$, then $G$ is Hamiltonian if $d(\nu) \ge n/2$ for all $\nu \in \mathcal{V}(G)$. 
\end{theorem}

The set of vertices $\mathcal{V}(G)$ is interpreted as the set of occupied orbital slots $\Tilde{\mathcal{J}} 	\subseteq \mathcal{J}$, defined as $\Tilde{\mathcal{J}}:=\{j:x_j=1,j \in \mathcal{J}\}$. The number of vertices equals the number of occupied orbital slots, and the degree $d(\nu)$ of a vertex is interpreted as the number of links that an occupied orbital slot $j \in \Tilde{\mathcal{J}}$ establishes at time step $t$. Then, we formulate Dirac's theorem as follows. First, constraint~\eqref{isl:cardinality_0} requires at least three occupied orbital slots. Second, constraints~\eqref{isl:vert_degree_0} enforce that each occupied orbital slot $j \in \Tilde{\mathcal{J}}$ must have a number of ISLs equal to at least half of the total number of occupied orbital slots.
\begin{subequations}
    \begin{alignat}{2}
    &\sum_{j\in \Tilde{\mathcal{J}}} x_{j} \ge 3 \label{isl:cardinality_0}\\
    &\sum_{j\in \Tilde{\mathcal{J}}\setminus\{k\}} W_{tjk}x_{j} \ge \frac{1}{2}\sum_{j\in \Tilde{\mathcal{J}}}x_{j}, \quad & \forall t \in \mathcal{T}, \forall k \in \Tilde{\mathcal{J}} \label{isl:vert_degree_0}
\end{alignat}
\end{subequations}
Therefore, any constellation configuration whose orbital slots satisfy these constraints will have a Hamiltonian cycle. However, determining the constellation configuration first, and then assessing whether it satisfies Dirac's theorem, is not efficient. To this end, we propose a reformulation of these constraints so that they consider the entire set of orbital slots $\mathcal{J}$, and that they can be adopted at the time of constellation configuration design. First, constraints~\eqref{isl:cardinality_0} are extended to the entire set of orbital slots $\mathcal{J}$ and reformulated as constraints~\eqref{isl:cardinality}. Second, constraints~\eqref{isl:vert_degree_0} are replaced by constraints~\eqref{isl:vert_degree}, which add an additional term that makes Dirac's theorem constraints nonbinding for non-occupied orbital slots.
\begin{subequations}
    \begin{alignat}{2}
    &\sum_{j\in \mathcal{J}} x_{j} \ge 3 \label{isl:cardinality}\\
        &J(1-x_k)+\sum_{j\in \mathcal{J}\setminus\{k\}} W_{tjk}x_{j} \ge \frac{1}{2}\sum_{j\in \mathcal{J}}x_{j}, \quad & \forall t \in \mathcal{T}, \forall k \in \mathcal{J} \label{isl:vert_degree}
\end{alignat}
\end{subequations}

\begin{remark}[One-Satellite Fault-Tolerant ISL Network Topology]
    The existence of a Hamiltonian cycle for each time step $t$ ensures that the ISL connectivity is robust against any one satellite failure. The reason is that removing a satellite from the cycle will yield a Hamiltonian path, that is, a path that visits all the satellites (vertices of the graph).
\end{remark}

\begin{remark}[One-Satellite Fault-Tolerant Network Design]\label{remark:robust_network}
    Suppose the mission designer aims to design a one-satellite fault-tolerant robust network, both in terms of ISL and satellite-target connectivity. In that case, we propose to adopt objective function~\eqref{sclp:obj} to minimize the constellation configuration cost, visibility constraints~\eqref{isl:visibility_robust} to enforce coverage requirements, and Dirac's theorem constraints~\eqref{isl:cardinality} and~\eqref{isl:vert_degree} to require the ISLs. Then, the formulation that determines the minimum-cost satellite network is given as:
    \begin{alignat}{2}
    \min \quad & \sum_{j\in \mathcal{J}}c_{j} x_{j} \notag \\
    \text{s.t.} \quad & \sum_{j\in \mathcal{J}} x_{j} \ge 3 \notag \\
            &J(1-x_k)+\sum_{j\in \mathcal{J}\setminus\{k\}} W_{tjk}x_{j} \ge \frac{1}{2}\sum_{j\in \mathcal{J}}x_{j}, &  \quad \forall t \in \mathcal{T}, \forall k \in \mathcal{J}  \notag \\
    & \sum_{j\in \mathcal{J}} V_{tjp}x_{j} \ge 2, &\quad  \forall t \in \mathcal{T}, \ \forall p \in \mathcal{P} \label{isl:visibility_robust}\\
    &x_{j} \in \{0, 1\}, &\quad \forall j\in \mathcal{J} \notag
    \end{alignat}
    The resulting satellite network is robust in terms of satellite-target and ISL connectivity. Forcing the right-hand side of constraints~\eqref{isl:visibility_robust} equal to two for all $t\in\mathcal{T}$ and for all $p\in\mathcal{P}$ makes the network robust in terms of satellite-target connectivity, as at least two satellites must establish links for each time step $t$. In addition, and as outlined above, the network is at least one-satellite fault-tolerant due to the existence of a Hamiltonian cycle. Lastly, this formulation extension can be adopted as a discipline for designing the constellation configuration in a general MDO architecture to tackle the constellation design. Leveraging this approach, the critical mission requirements that ensure the feasibility of the ISL topology, such as considering the Doppler shift, link budget, transmission latency thresholds, maximum power consumption, and data downlink requirements, would be addressed by the telecommunications discipline; thereby, both disciplines would closely collaborate under the general coordination of the MDO system optimizer.
\end{remark}

\begin{example}[One-Satellite Fault-Tolerant Network Design]
To showcase Remark~\ref{remark:robust_network}, we propose an illustrative example that designs a robust satellite constellation for uninterrupted communications between Morgantown, West Virginia, and Brisbane, Queensland, Australia (geodetic coordinates \{$27.46\degree$S, $153.02\degree$E, \SI{0}{m}\}). For both targets, the minimum elevation angle is set to \SI{5}{deg}. The epoch is January 1, 2025, 12:00:00.000 UTC, and the time step size is \SI{240}{s}, selected as a trade-off between the problem scalability and modeling fidelity. All orbital slots belong to circular RGT orbits with a 7:1 resonance ratio, a repetition period of \SI{85951.43}{s}, a semi-major axis of \SI{11507.30}{km}, and an inclination of \SI{45}{deg}.

The optimal constellation configuration that delivers continuous two-fold coverage over the ground targets is obtained after \SI{885.65}{s}, and uses 13 satellites. The leftmost column of Fig.~\ref{fig:sclp_isl_2d} outlines the constellation pattern vector $\bm{x}$ and the spatial distribution of the orbital slots in the RAAN versus argument of latitude plane. Further, the rightmost column of Fig.~\ref{fig:sclp_isl_2d} shows the constellations' reference visibility profile $\bm{v}_{p}$ and coverage timeline $\bm{b}_{p}$. Lastly, Fig.~\ref{fig:sclp_isl_3d} illustrates the 3D visualization of the optimal constellation configuration, ISLs, and targets at the epoch.

\begin{figure}[htpb]
    \centering
    \includegraphics[width=0.9\linewidth]{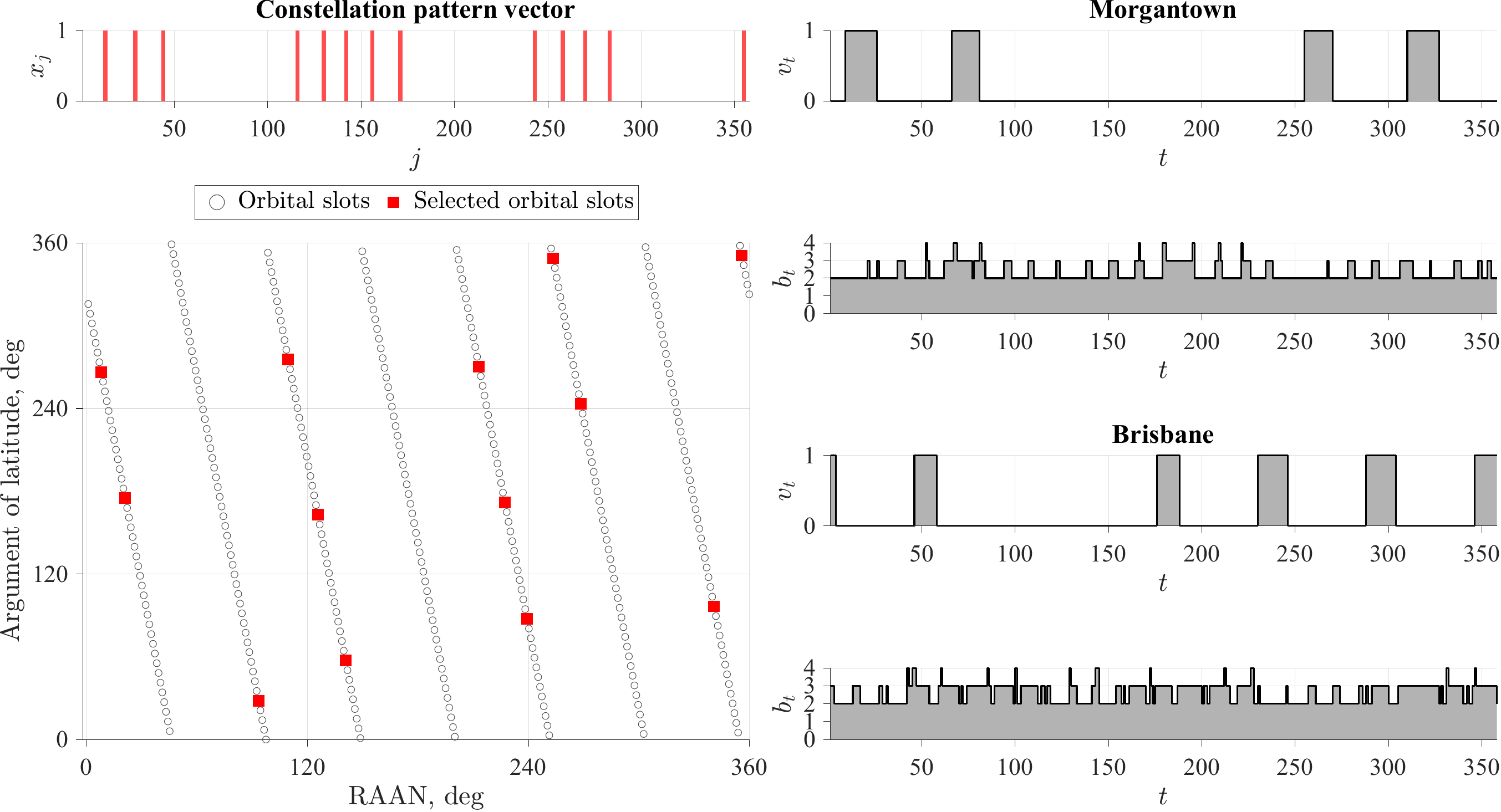}
    \caption{APC decomposition and orbital slots' distribution.}
	\label{fig:sclp_isl_2d}
\end{figure}

\begin{figure}[htpb]
    \centering
    \includegraphics[width=0.8\linewidth]{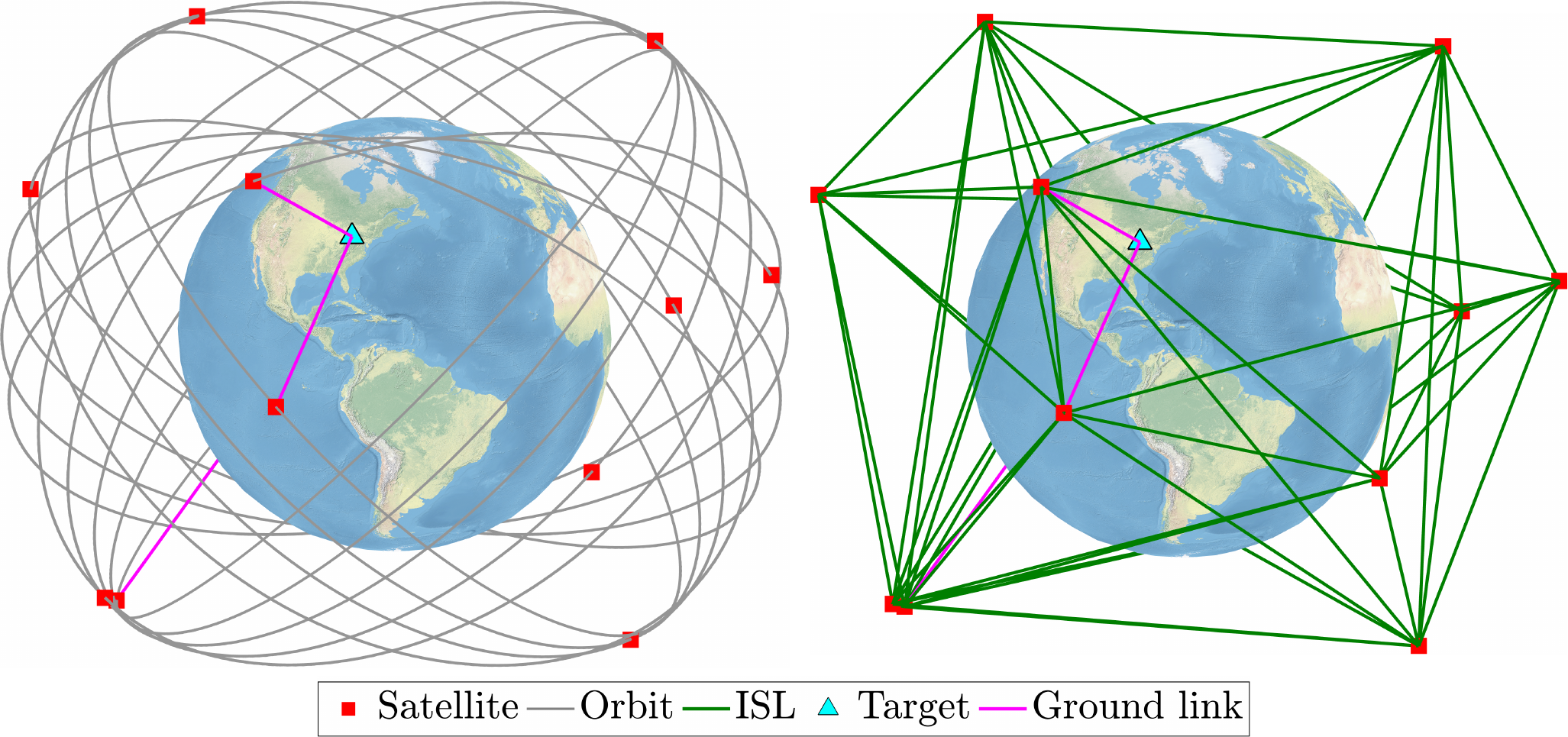}
    \caption{Thirteen-satellite optimal constellation configuration.}   
    \label{fig:sclp_isl_3d}
\end{figure}

To illustrate the robustness of the constellation against a one-satellite failure, we randomly remove one. Figure~\ref{fig:sclp_isl_pat_cov_comp} presents the coverage timelines over Morgantown and Brisbane after removing the satellite located in orbital slot $j=243$. The coverage drops from a minimum of double-fold to single-fold at specific time steps due to the failure. Further, owing to Dirac's theorem, the constellation remains fully connected after the failure. Figure~\ref{fig:sclp_isl_failure} presents the network topology for two time steps, $t=119$ and $t=257$, projected over the latitude-longitude plane. The first column of the figure enables the visualization of the impact on the number of connections by randomly removing a satellite, leaving Brisbane with a single link to a satellite at time step $t=119$. Similarly, the second column of the figure depicts the reduction in ground-to-satellite links suffered by Morgantown at time step $t=257$.

\begin{figure}[htbp]
    \centering
    \includegraphics[width=0.492\linewidth]{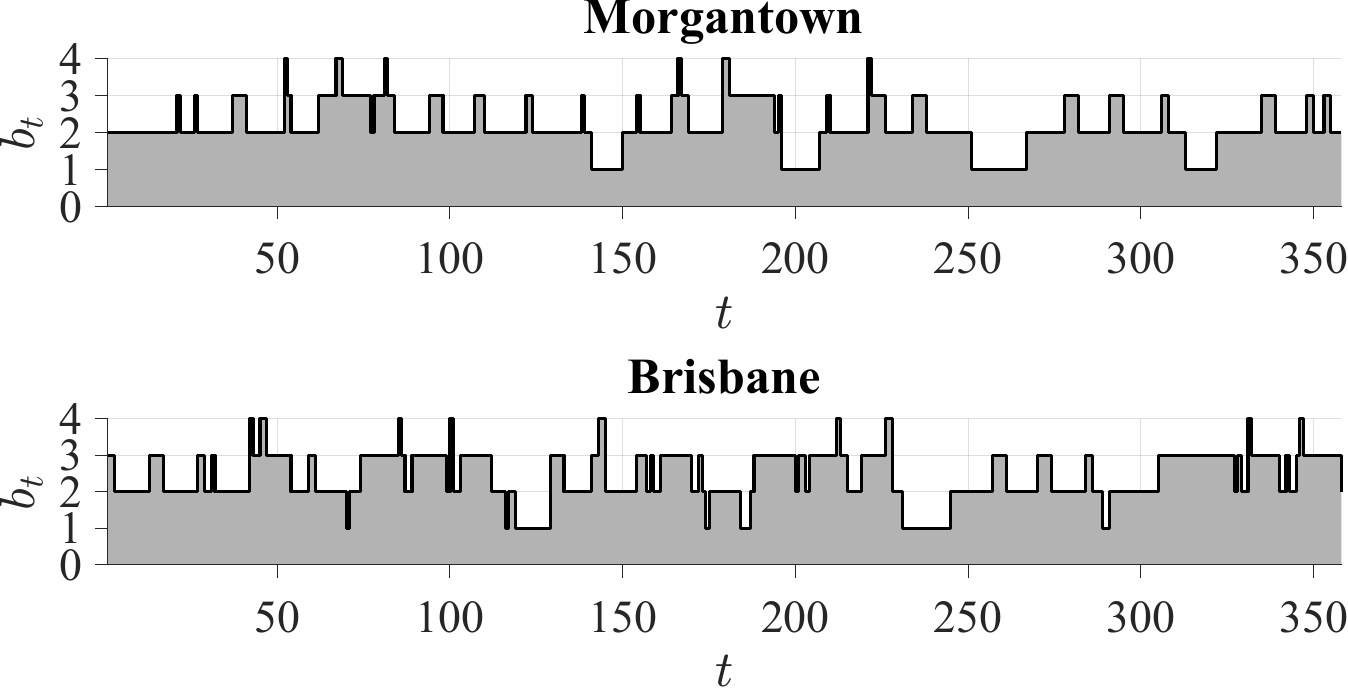}
    \caption{Coverage timelines after one-satellite failure.}
    \label{fig:sclp_isl_pat_cov_comp}	
\end{figure}

\begin{figure}[htbp]
    \centering
    \includegraphics[width=0.82\linewidth]{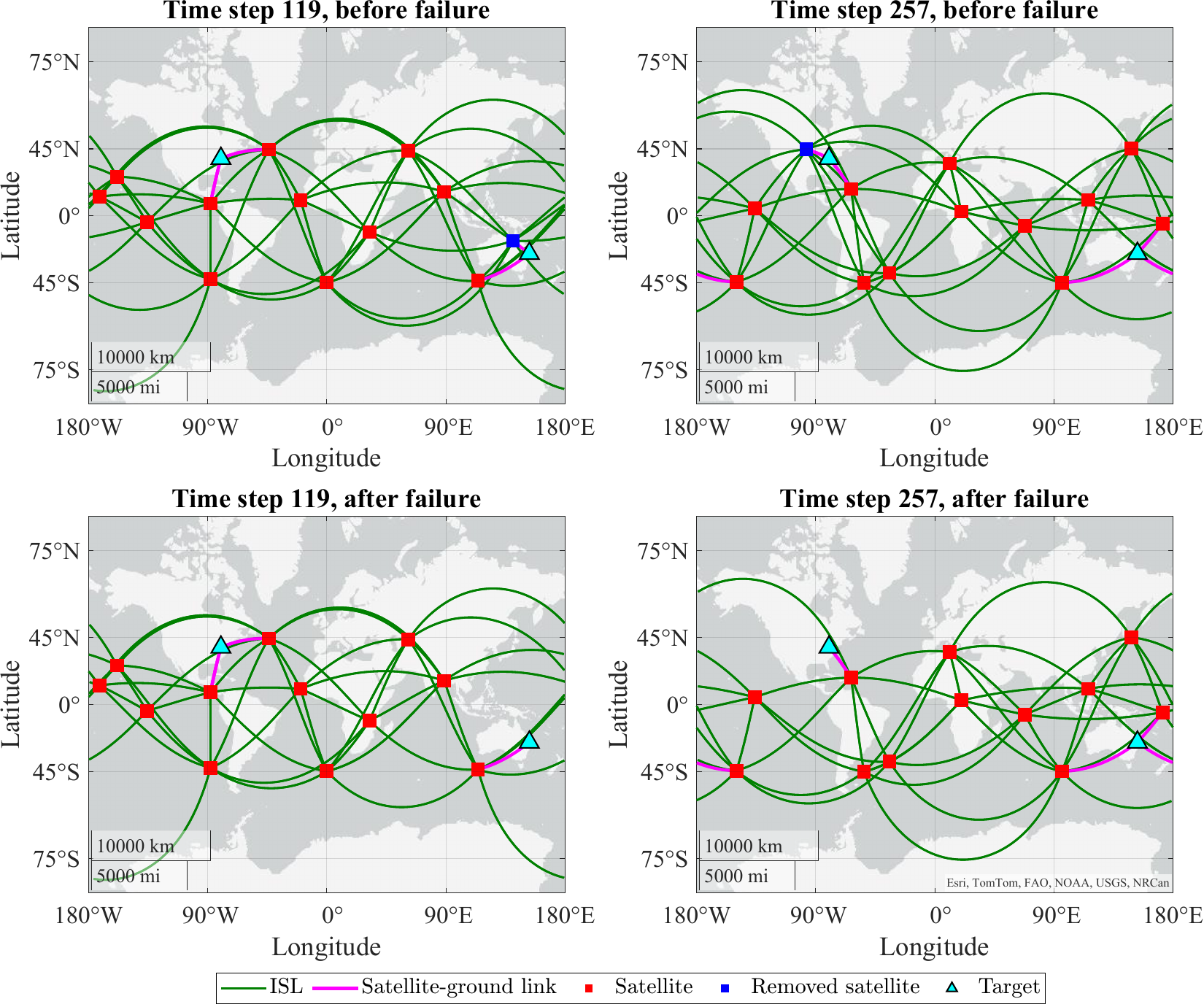}
    \caption{Projected view of optimal constellation configuration, and ISLs before and after one-satellite failure.}
    \label{fig:sclp_isl_failure}
\end{figure}
\end{example}

\clearpage
\section{Conclusions}\label{sec:conclusions}

This paper tackles the problem of fragmented and problem-specific constellation configuration design methodologies offered in the literature by proposing a unified and coherent collection of MILP. Each formulation determines a provably-optimal constellation configuration adopting key figures of merit (\textit{e.g.,} constellation configuration cost, percentage coverage, MRT, ART). These figures of merit act as the optimization objective, or as constraints based on user preferences, thereby configuring the MILP formulation scope. In addition, we extend the scope of the collection by deriving new MILP formulations suitable for initially overlooked mission scenarios.

To showcase the applicability of the proposed formulations, we present a set of illustrative examples, comparative analyses, mix-and-match operations, and formulation extensions. First, assuming a common set of parameters, the illustrative examples elucidate the objective of each formulation, which are then contrasted in terms of the required number of satellites, percentage coverage, MRT, and ART throughout a set of comparative analyses. Second, the mix-and-match operations interchange objective function and constraints from the MILP collection to accommodate new mission requirements (\textit{e.g.,} MRT or ART as mission requirements). Third, the formulation extensions introduce novel constraints or interpret differently the set of parameters to further address new mission specifications (\textit{e.g.,} ISL, SSA constellation design with dynamic targets).

The contributions presented set the stage for future research endeavors. First, the formulations presented are not meant to be comprehensive; therefore, the mix and match, as well as the presented extensions, can be further expanded to enable mission designers to incorporate additional mission objectives or requirements that can be linearly represented into the constellation configuration optimization (\textit{e.g.,} orbital slot-dependent observation reward, alternative approaches to model the sum of MRTs such as harmonic mean). Second, alternative parameter generation methods can be investigated to optimize the computational runtime dedicated to their generation. Third, considering that all formulations are MILP, computational tractability analyses could be performed to characterize the scalability of the formulations relative to instance size (\textit{e.g.}, parameter dimensionality) and explore solution methods to enhance the efficiency of the proposed approaches. Lastly, new modeling approaches could be explored and compared to the ones proposed in this paper (\textit{e.g.,} formulations with tighter Big-M method constants).

\appendix
\section{Repeating Ground Track Orbital Slots Generation}\label{app:oe_generation}

RGT orbits are characterized by a repetition period that depends on the ratio between the satellite's number of ascending node crossings $N_{\text{P}}$ in $N_{\text{D}}$ nodal periods of Greenwich. Formally, the repetition period $T_{\text{R}}$ is defined as:
\begin{equation}
T_{\text{R}}=N_{\text{P}}T_{\text{S}}=N_{\text{D}}T_{\text{G}}
\end{equation}
where $T_{\text{S}}$ and $T_{\text{G}}$ are the satellite and Greenwich's nodal periods, respectively, defined as:
\begin{subequations}
    \begin{alignat}{2}
        &T_{\text{S}} = \frac{2\pi}{\dot{\omega} + \dot{M}}\\
        &T_{\text{G}} = \frac{2\pi}{\omega_{\oplus} - \dot{\Omega}}
    \end{alignat}
\end{subequations}
where $\dot{\omega}$ is the argument of perigee's rate of change due to perturbations, $\dot{M}$ is the mean anomaly's rate of change due to perturbations and nominal motion, $\omega_{\oplus}$ is the Earth's rotation rate, and $\dot{\Omega}$ the orbit's nodal regression rate of change. These parameters are defined as:
\begin{subequations}
    \begin{alignat}{2}
        & \dot{\omega} = \frac{3}{2}J_{2}\left( \frac{R_{\oplus}}{a(1-e^{2})} \right)^{2} \left(\frac{\mu_{\oplus}}{a^{3}}\right)^{1/2} \left[ 2-\frac{5}{2}\sin^{2}{i} \right]\\
       & \dot{M} = \left(\frac{\mu_{\oplus}}{a^{3}}\right)^{1/2} \left[ 1-\frac{3}{2}J_{2}\left( \frac{R_{\oplus}}{a(1-e^{2})} \right)^{2} (1-e^{2})^{1/2}\left( \frac{3}{2}\sin^{2}{i} -1 \right)\right]\\
       &\dot{\Omega} = -\frac{3}{2}J_{2}\left( \frac{R_{\oplus}}{a(1-e^{2})} \right)^{2} \left(\frac{\mu_{\oplus}}{a^{3}}\right)^{1/2} \cos{i}
    \end{alignat}
\end{subequations}
where $\mu_{\oplus}$ is the Earth's standard gravitational parameter, and $a$, $e$, and $i$ are the orbit's semi-major axis, eccentricity, and inclination, respectively. Then, given a user-defined resonance ratio $N_{\text{P}}:N_{\text{D}}$, inclination $i$, and eccentricity $e$, the semi-major axis $a$ is obtained by applying the Newton-Raphson method proposed in Ref.~\cite{bruccoleri_phd}.

Given that all orbital slots $j \in \mathcal{J}$ must share the same semi-major axis $a$, eccentricity $e$, and inclination $i$, the only degrees of freedom for orbital slot $j$ are the RAAN $\Omega_{j}$ and mean anomaly $M_{j}$, which can be obtained by satisfying the following condition \cite{lee2020satellite}:
\begin{equation}
    N_{\text{P}}\Omega_{j}+N_{\text{D}}M_{j}=\text{constant mod}(2\pi)
\end{equation}

\section*{Acknowledgment}
This research is based on work sponsored by TelePIX Co., Ltd. The authors would like to thank the anonymous reviewers whose comments increased the clarity and quality of this paper.

\bibliography{references}

\end{document}